\documentclass[twoside]{article}

\usepackage[accepted]{aistats2022}
\newlength\figH
\newlength\figW
\usepackage{graphicx}
\usepackage{pgfplots}
\usepgfplotslibrary{groupplots}
\usepackage{wrapfig}
\usepackage{nameref}
\usepackage{blindtext}
\usepackage{caption}
\usepackage{subcaption}

\usepackage{afterpage}

\usepackage{hyperref}
\hypersetup{
    colorlinks = true,
}

\usepackage[capitalise]{cleveref}

\usepackage{algorithm}
\usepackage{algorithmic}
\usepackage{epigraph}
\usepackage{tikz}
\usepackage{amsthm}
\usepackage{amsfonts}
\usepackage{amssymb}
\usepackage{graphicx}
\usepackage[all]{xy}
\usepackage{amsmath}
\usepackage{mathtools}
\usepackage{pgfplots}
\pgfplotsset{compat=1.11}
\usepackage{setspace} 
\usepackage{wasysym}
\usepackage{cite}
\usepackage{xfrac}
\usepackage{float}
\usepackage{soul,xcolor}
\setstcolor{red}
\usetikzlibrary{backgrounds,automata}

\usetikzlibrary{arrows,shapes,positioning,math,patterns}
\usetikzlibrary{decorations.markings}
\usetikzlibrary{calc}
\tikzstyle arrowstyle=[scale=1]
\tikzstyle directed=[postaction={decorate,decoration={markings, mark=at position .65 with {\arrow[arrowstyle]{stealth}}}}]

\DeclareMathOperator{\po}{po}
\DeclareMathOperator{\cho}{cho}

\DeclareMathOperator{\conco}{con}
\DeclareMathOperator{\disco}{dis}

\DeclareMathOperator{\divers}{div}

\DeclareMathOperator{\myand}{and}

\newtheorem{mydef}{Definition}

\newtheorem{mypro}{Proposition} 
\newtheorem*{myfpr}{Maximum Choice Problem} 
\newtheorem*{mychar}{Characterization of Hyperbolas} 
\newtheorem*{myifo}{Integral Formula} 
\newtheorem{myfle}{Lemma}

\newtheorem{mytheo}{Theorem} 
\newtheorem{mycoro}{Corollary}

\begin{document}

\afterpage{\cfoot{\thepage}}
\clearpage

%

%

\twocolumn[

\aistatstitle{$k$-Pareto Optimality-Based Sorting with Maximization of Choice}

\aistatsauthor{ Jean Ruppert* \And Marharyta Aleksandrova** \And  Thomas Engel** }

\aistatsaddress{ \scriptsize{*Mathematics and Computing S.à.r.l., \url{first.last@mathcomp.lu}} \And  \scriptsize{**University of Luxembourg, \url{first.last@uni.lu}} 
} ]




\begin{abstract}

Topological sorting\footnote{
\textit{Topological sorting}\cite{Taocp} here means the process of sorting a set of items with respect to a preference or dominance relation. We use the terms \textit{preference} and \textit{dominance} interchangeably.
}
is an important technique in numerous practical applications, such as information retrieval, recommender systems, optimization, etc. 
In this paper, we introduce a problem of \textit{generalized topological sorting with maximization of choice}, that is, of choosing a subset of items of a predefined size that contains the maximum number of equally preferable options (items) with respect to a dominance relation. 
We formulate this problem in a very abstract form and prove that sorting by \textit{$k$-Pareto optimality} yields a valid solution. 
Next, we show that the proposed theory can be useful in practice. We apply it during the selection step of genetic optimization and demonstrate that the resulting algorithm outperforms existing state-of-the-art approaches such as \textit{NSGA-II} and \textit{NSGA-III}. 
We also demonstrate that the provided general formulation allows discovering interesting relationships and applying the developed theory to different applications.

\end{abstract}

\section{Introduction}
\label{sec:introduction}

In the modern era of information overload, the task of 
choosing a subset of the most useful items 
is extremely important. 
Various tools were developed with the aim to assist a user with this task, for example, text search engines \cite{croft2010search} and recommender systems \cite{resnick1997recommender}.
In most of the cases, such systems suggest to the user a small set of elements\footnote{
We use the terms \textit{item} and \textit{element} interchangeably.
}.
Thereby, if the number of equally preferable options is large, a heuristic is used to discard a fraction of them.
However, in some applications the user might be willing to analyze all equally preferable options with the aim to choose the best one.
This can happen, for example, in the case of choosing a habitation.

A similar problem of \textit{choosing a subset of most preferable elements} also arises as an important step when solving various practical tasks.
A straightforward example would be the selection step in genetic optimization algorithms~\cite{mitchell1998introduction}.
At this step, a subset of the current population is chosen to advance to the next generation.
Having the chosen subset made up of elements with large fitness values guides the evolution process in the desired direction. 
At the same time, selecting a subset with the largest variety of genes ensures variability of characteristics and allows faster exploration of the search space.

These examples bring us to the problem of \textit{generalized\footnote{
Later we show that the formulated problem is a generalization of topological sorting. In the case of a partial order relation, standard topological sorting is considered.
} topological sorting with choice maximisation} which we also refer to as \textit{maximum choice problem}.
This problem aims to choose a subset of a predefined maximum size,  consisting of   most preferable items and containing the largest number of equally preferred elements\footnote{
If items are equally preferable, then we say that they offer choice for the user or the system.
}. 
To the best of our knowledge, this problems has not yet been studied in the literature.
In this text, we propose a theoretical solution to the maximum choice problem and demonstrate how both the problem and its solution can be applied in practice.

The contributions of this work are the following: 
\begin{enumerate}
    \item We formulate the maximum choice problem in a broad sense for arbitrary elements, preference
    relations $R$, and measures $\mu$ indicating the set sizes,
    see \cref{sec:pareto-optimality}.
    
    \item We propose a solution  based on the concept of $k$-Pareto optimality, whose definition relies on the 
    relation $R$, see \cref{posolutions}.
    
    \item We further investigate the proposed solution from a theoretical point of view and discover interesting characteristics, such as the relationship between $k$-Pareto optimal elements and the arc of hyperbola, see \cref{sec:theoretical-exploration}.

    \item Finally, we demonstrate the applicability of our approach to real-world problems by considering genetic optimization,
    see \cref{sec:practical-exploration}.
\end{enumerate}

\section{Sorting with Choice Maximization}
\label{sec:pareto-optimality}

To formally define the problem of \textit{maximization of choice}, we introduce several definitions in \cref{sec:pareto-optimality:definitions,sec:pareto-optimality:max-choice}.
The resulting formalization of this problem is abstract and quite general.
However, this generality allows discovering novel connections and applying the developed theory to numerous practical problems.
We also illustrate the defined concepts with an example in \cref{sec:pareto-optimality:examples}.

\subsection{Definitions}
\label{sec:pareto-optimality:definitions}

We consider a set $X$ with 
a binary relation $R$. Intuitively $x R y$ means that $x$ is preferable to $y$. 
The case $xRy$ and not $yRx$ means that $x$ is \emph{strictly preferable} to $y$. This situation is denoted by $x R^* y$.
We also consider a positive and $\sigma$-finite measure $\mu$ \cite{halmos2013measure} defined on $X$.
Thus, we have a measure space $(X, \Sigma,\mu)$,  where $\Sigma$ is a set of subsets of $X$, and $\mu$ intuitively indicates the size of these subsets.  
To ensure measurability, throughout this text the  characteristic function $1_R$  of the relation $R$ is assumed to be sufficiently regular\footnote{
$1_{R}$ is equal to $1$ if $xRy$ and is equal to $0$ otherwise.}.
The $\mu$ can be defined in different ways. 
Important examples are the counting measure and  probability measures $P$.
Depending on the definition of $\mu$, it can indicate the following characteristics of the elements in $X$:
\textit{how many?}, \textit{how likely?}, \textit{how important?}, or \textit{what volume?}

To illustrate these definitions, we  consider the following example. 
Let $X$ be a set of possible habitations of which the user has to choose the best according to his
preferences encoded by the relation $R$. 
In such a situation, the relation $R$ can be multidimensional. 
Let us assume, for simplicity, that an optimal habitation for the user is close to a given location, for example, his workplace (relation $R_l$), is situated in a district with a smaller population size (relation $R_p$), and is close to a river (relation $R_r$).
Thus, the user's preferences can be represented by the preorder relation $R= R_l \& R_p \& R_r$.
In our example, all available habitations from $X$ can be mapped onto points in a 3-dimensional space of \textit{Proximity to the location} $\times$ \textit{Population} $\times$ \textit{River}.
The fact that $R$ is a preorder relation means that some elements of $X$ can be comparable, while others not.
For example, the habitation $x$ with coordinates $(50, 100, True)$ is strictly preferable to $y$ with coordinates $(60, 100, True)$, that is $x R^* y$.
At the same time, the habitation $z$ with coordinates $(40, 100, False)$ is incomparable with $x$.
Indeed, $z$ is better with respect to $R_l$, it is situated closer to the required location, but $x$ is better with respect to $R_r$, as the latter is situated near a river.

Having the task to find a subset of $X$ that is `best' according to $R$, a rational solution can be formulated with the following recursive expression: \textit{if an element $x$ is selected, then all elements that are strictly preferable to $x$ should be also selected}.
In our example, this translates into the task of finding a subset of habitations $S_R$ that might be suitable for the user.
Naturally, if $y \in S_R$, then $x \in S_R$ as the latter corresponds better to the preferences of the user defined by the relation $R$.
We formalize this rationality condition by defining \textit{selections} as follows.

\begin{mydef}
\label{def-selection}
	A \emph{selection} $S$ is a subset of $X$ such that $x \in S$ and $y R^* x$ implies $y \in S$. The set of all selections in $\Sigma$ is denoted by $\mathcal{S}$\footnote{
	If $R$ is a partial order relation, then the selections are the down-sets \cite{davey2002introduction}.
	}.
\end{mydef}

\subsection{The Maximum Choice Problem}
\label{sec:pareto-optimality:max-choice}

As discussed in \cref{sec:introduction}, in practical applications when selecting a subset of $X$ one might want not only to respect the above rationality constraint, but also to maximize the number of incomparable pairs.
The latter condition is equivalent to the maximization of the diversity of the selected subset,
or the maximization of the provided choice. 
In our example with habitations, if both $x$ and $z$ are presented to the user, then he can
choose an appropriate habitation by himself\footnote{
We consider the case when addition preferences cannot be encoded and the user has to make the final choice.
}.

In terms of our notations, this will be translated into the condition of 
selecting as many pairs $x, y$ 
such that neither $x$ is strictly preferable to $y$ ($ \neg x R^* y$) nor $y$ is strictly preferable to $x$ ($\neg y R^* x$). 
This means that there is freedom of choice between $x$ and $y$ ($xRy=yRx$). 
This motivates the following quantitative definition of choice for  measurable subsets of $X$.
\begin{mydef}
\label{def:choice}
		\emph{Choice} offered by a 
		set $A$ is the number
    	$\cho(A)=(\mu \times \mu)(\{(x, y) \in A^2 | xRy=yRx\})$.
\end{mydef}
The choice offered by a measurable set $A$ essentially measures how many pairs of items offering choice can be extracted from $A$.
Additionally, if one wants to restrict the size of the selected subset, in our example, to present to the user a small set of suitable habitations with $\mu(S_R) \leq m$,
then this leads us to the definition of the maximum choice problem:
\begin{myfpr}
	For a given $m$ find all selections $T$ such that
	$\cho(T) 
	=\max\limits_{\substack{ S \in \mathcal{S}, \mu(S) \leq m}} \cho(S).$
	Any such selection $T$ will be said to \emph{offer maximum choice for $m$}.
\end{myfpr}

In practical applications, it might be more insightful to consider the concept of \textit{diversity} that is functionally related to the concept of choice.
\begin{mydef}
\label{def:diversity}
For any measurable set $A$, the \em{diversity} of $A$ is the ratio
$\divers(A)=\cho (A) / \mu(A)^2.$
\end{mydef}  
Thus $\divers(A)$ is the likelihood that there is choice between the two randomly chosen elements. 
The main ingredient of our solution to the maximum choice problem is the following concept.
\begin{mydef}
\label{def-po}
	The \emph{$k$-Pareto optimality} of an element $x \in X$ is the measure of the subset of $X$ containing all elements strictly preferable to $x$:
	$\po(x)=\mu(\{y| yR^*x\})$.
\end{mydef} 
If $\mu$ is the probability measure, the $k$-Pareto optimality of an element $x$, $\po(x)$, is the likelihood an element drawn at random from $X$ is strictly preferable to $x$.

Finally, we introduce the sets of \textit{at least $k$-Pareto optimal elements}. Such sets consist of all elements in $X$ with $\po \leq k$.
In \cref{posolutions} we show that this concept yields a solution for 
the maximum choice problem.  
\begin{mydef}
\label{def:po-set}
	The \emph{at least $k$-Pareto optimal elements} $T_k$ form the measurable set defined as follows:\\
	$T_k=\{x \in X | \po(x) \leq k  \}$\footnote{
    If $R^*$ is transitive, then for any $k$, $T_k$ is a selection.	
	}.
\end{mydef}

\subsection{Example}
\label{sec:pareto-optimality:examples}

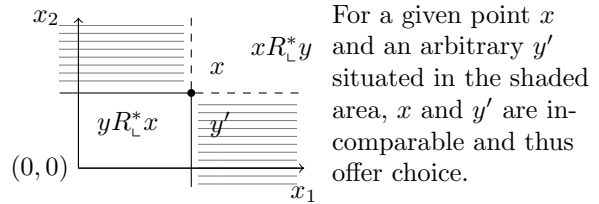
\begin{figure}
  \centering
  \tikzmath{\ssw = 3;  
\ssh = 2;  
 \x1 =1.5; \x2=\x1 - 0.5; \fillepsilon=0.1; 
 \excess=0.25; } 

\begin{tikzpicture}
\coordinate (x) at(\x1,\x2);
\coordinate (rect11) at (\ssw,\ssh); 
\coordinate (rect10) at (\ssw,0); 
\coordinate (rect01) at (0,\ssh); 
\coordinate (rect00) at (0,0); 
\coordinate (leftDown) at(0.3,-0.2);

\coordinate (txtChoice) at ($(rect11)!0.5!(rect10)+(\excess, 0)$);

\fill[pattern=horizontal lines, pattern color=gray!100]       
($(0,\ssh)+(-\excess,-\fillepsilon)$) 
-- ($(0,\x2)+(-\excess,\fillepsilon)$) 
-- ($(x)+(-\fillepsilon, \fillepsilon)$)
--($(\x1,\ssh)+(-\fillepsilon,-\fillepsilon)$) ;

\fill[pattern=horizontal lines, pattern color=gray!100]       
($(x)+(\fillepsilon,-\fillepsilon)$) 
-- ($(\x1,0)+(\fillepsilon,-\excess)$) 
-- ($(rect10)+(-\fillepsilon, -\excess)$)
--($(\ssw,\x2)+(-\fillepsilon,-\fillepsilon)$) ;

\draw (rect00) to (rect10);
\draw (rect00) to  (rect01);

\draw  [->](0,0) coordinate (Oo) -- (\ssw, 0) coordinate (Ox) ;
\draw  [->](Oo) -- (0, \ssh) coordinate (Oy) ;
\draw ($(0,\x2)-(\excess, 0)$) to (x) to ($(\x1,0)-(0,\excess)$);
\draw [dashed](\x1,\ssh) to (x) to (\ssw,\x2);

\fill (x) circle [radius=1.5pt];

\node [label=above right:{$x$}] at (x) {};
\node [label=below right:{$y R^{*}_{\llcorner} x$}] at (0,\x2) {};
\node [label=below:{$x R^{*}_{\llcorner} y$}] at ($(rect11)+(-0.3,0)$) {};
\node[anchor=west,text width=3.5cm] at (txtChoice) {
For a given point $x$ 
and an arbitrary $y'$ 
situated in the shaded  
area, $x$ and $y'$ are 
incomparable and thus 
offer choice.};
\node [label=below right:{$y'$}] at (x) {};

\node [left] at (rect00){$(0,0)$};
\node [label=left:{$x_2$}] at (rect01) {};
\node [label=below:{$x_1$}] at (rect10) {};
\end{tikzpicture}
  \caption{Illustration of the relation $R_{\llcorner}$.}
  \label{fig:R_llcorner}%
\end{figure}

In this subsection, we discuss an illustrative example to demonstrate the concepts defined in \cref{sec:pareto-optimality:definitions,sec:pareto-optimality:max-choice}.
Let us consider a finite subset $X$ of $\mathbb{R}^2$, counting measure $\mu$, and the relation $R_{\llcorner}$ defined as follows: $(x_1, x_2) R_{\llcorner} (y_1, y_2)$ iff $x_1 \leq y_1$  and  $x_2 \leq y_2$, see \cref{fig:R_llcorner}.
In economics, $x$ is Pareto optimal if there is no $y$ in $X$ such that $x R_{\llcorner}^* y$.
In our language, this means that $\po(x)=0$. 
Thus, $k$-Pareto optimality indicates how much an element is away from being Pareto optimal.

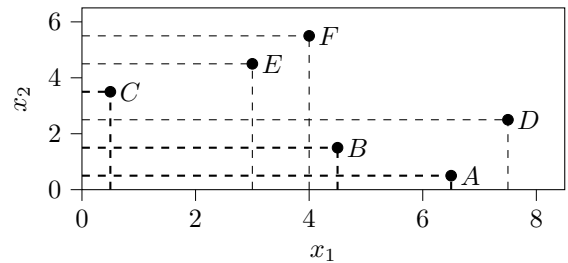
\begin{figure}
	\centering{
\begin{tikzpicture}

\begin{axis}[
tick align=outside,
tick pos=left,
x grid style={white!69.0196078431373!black},
xmin=0, xmax=8.5,
xtick style={color=black},
y grid style={white!69.0196078431373!black},
ymin=0, ymax=6.5,
ytick style={color=black},
height=4.cm,
width=8.cm,
xlabel={$x_1$},
ylabel={$x_2$},
]

\filldraw[black] (6.50,0.50) circle (2pt)
node[anchor=west] {$A$};
\filldraw[black] (4.50,1.50) circle (2pt)
node[anchor=west] {$B$};
\filldraw[black] (0.50,3.50) circle (2pt)
node[anchor=west] {$C$};
\filldraw[black] (7.5,2.50) circle (2pt)
node[anchor=west] {$D$};
\filldraw[black] (3.0,4.50) circle (2pt)
node[anchor=west] {$E$};
\filldraw[black] (4.0,5.50) circle (2pt)
node[anchor=west] {$F$};

\addplot [draw=black, fill=black, mark=., only marks, scatter]
table{%
x  y
6.50 0.50 0
4.50 1.50 0
0.50 3.50 0
7.50 2.50 0
3.0 4.50 0
4.50 5.50 0
};

\addplot [thin, black, dashed]
table {%
0 2.50
7.50 2.50
};
\addplot [thin, black, dashed]
table {%
7.50 0
7.50 2.50
};
\addplot [thick, black, dashed]
table {%
0 0.50
6.50 0.50
};
\addplot [thick, black, dashed]
table {%
6.50 0
6.50 0.50
};
\addplot [thick, black, dashed]
table {%
0 1.50
4.50 1.50
};
\addplot [thick, black, dashed]
table {%
4.50 0
4.50 1.50
};
\addplot [thin, black, dashed]
table {%
0 5.50
4.0 5.50
};
\addplot [thin, black, dashed]
table {%
4.0 0
4.0 5.50
};
\addplot [thin, black, dashed]
table {%
0 4.50
3.0 4.50
};
\addplot [thin, black, dashed]
table {%
3.0 0
3.0 4.50
};
\addplot [thick, black, dashed]
table {%
0 3.50
0.50 3.50
};
\addplot [thick, black, dashed]
table {%
0.50 0
0.50 3.50
};

\end{axis}

\end{tikzpicture}
		}
    \caption{Computation of Pareto-optimality. Example of a finite subset of $\mathbb{R}^2$, the counting measure, and $R_{\llcorner}$. 
    }
    \label{fig:po-example} 
\end{figure}

Let $X$ be comprised of six points presented in \cref{fig:po-example}.
As we are considering the counting measure, $\mu (X) = 6$.
Points $A$, $B$, and $C$ are not dominated by any other point.
This means that $\po(A) = \mu(\{y| yR^*A\}) = 0$ and $\po(B) = \po(C) = \po(A) = 0$.
Point $E$ is dominated by a single point $C$, that is $\po(E) = \mu(\{C\}) = 1$.
Finally, points $F$ and $D$ are dominated by two other points each, resulting in $\po(F) = \mu(\{ C, E \}) = 2$ and $\po(D) = \mu(\{ A, B \}) = 2$.

Sorting the set $X$ by $k$-Pareto optimality of its elements will produce the following result: $(\{A, B, C\}, \{E\}, \{D, F\})$.
This sorting is different from sorting by Pareto fronts.
The latter approach is widely used in practice and is the basis of all Pareto dominance-based genetic optimization algorithms \cite{li2015many}.
Sorting by Pareto fronts is done in the following way.
First, the first Pareto front, which is the set of non-dominated points, is identified.
Next, the points from this front are removed from the consideration and the process is repeated until all points are assigned to a front.
Sorting the points from \cref{fig:po-example} by Pareto fronts will produce the following result: $(\{A,B,C\}, \{D, E\}, \{F\})$. 
Note, that the point $D$ moved from the 3d equivalence class when sorting by $k$-Pareto optimality to the 2d when sorting by Pareto fronts.
Let us consider the two selections of size 4: $S_E = \{A,B,C, E\}$ and $S_D = \{A,B,C, D\}$.
Sorting by Pareto fronts does not distinguish between these two selections as both $E$ and $D$ belong to the same equivalence class.
At the same time, sorting by $k$-Pareto optimality has a larger preference towards $S_E$. 
Also, the latter selection contains more incomparable pairs of elements and thus offers more choice.
Indeed, $D$ is incomparable with only one point $C$, but $E$ is incomparable with two points: $A$ and $B$.

\section{The Maximum Choice Theorem}
\label{posolutions}

The main result of this text states that at least $k$-Pareto optimal elements ($T_k$) are the largest measurable selections offering maximum choice for their respective measures.

\begin{mytheo}
\label{theorem:po-solution}
A set of at least $k$-Pareto optimal elements $T_k$ offers maximum choice for $\mu(T_k)$ if it is a selection and if $\mu(T_k)<+\infty$.
Moreover $T_k$ is the largest selection offering maximum choice for $\mu(T_k)$ in the sense that it contains any other selection offering maximum choice for $\mu(T_k)$. 
\end{mytheo}

The second part of the above theorem  precisely means that if a selection $A$ offers maximum choice for $\mu(T_k)<+\infty$, then $A \subseteq T_k$ almost-everywhere. One may note that the above theorem does not state that the at least $k$-Pareto optimal elements are the only largest selections offering maximum choice. Later, we give an example of a fundamentally different selection offering maximum choice for a value of $m$ where there is no such $k$ that $\mu(T_k)=m$. 

We prove the above stated theorem in several steps. 
First, we show that for selections the computation of choice can be simplified. It only requires to compute a simple integral instead of a double integral.
\begin{myifo} 
If $S$ is a measurable selection and $\mu(S)<+\infty$, then 
\begin{align} \label{iFo}
\cho(S) & = \int_{S} (\mu(S)-2\po(x)) d\mu(x). 
\end{align}
\end{myifo}

\begin{proof}
The fact that set $S$ 
is a selection means that $\forall y \in S:$ $\{x \in S | x R^* y\}=\{x \in X | x R^* y\}$.
That is, any element $x$ from $X$ strictly preferable to any element $y$ in $S$, also belong to $S$ ($x \in S$).
Using the definition of choice from Def.~\ref{def:choice} and $\mu(S)<+\infty$ we obtain 
\begin{align*}
\cho(S) &=\mu(S)^2 - 2(\mu \times \mu)(\{(x,y)\in S^2 | x R^* y\}) \\
  &=\mu(S)^2 - 2 (\mu \times \mu)(\{(x,y)\in S \times X | y R^* x\}).
\end{align*}
Fubini's theorem \cite{halmos2013measure} indicates that
\begin{equation*}
\begin{split}
 (\mu \times & \mu)(\{(x,y)\in S \times X | y R^* x\} = \\
 &= \int_{S}\left ( \int_{X} 1_{R^*} d(\mu(y)) \right )d\mu(x) = \int_{S} \po(x) d\mu(x), 
\end{split}
\end{equation*}
where $1_{R^*}$ is the characteristic function of $R^*$\footnote{
$1_{R^*}$ is defined on $X^2$ in a similar way to $1_{R}$: it is equal to $1$ if $xR^*y$ and is equal to $0$ otherwise.
}.

Finally, the integral formula results from the fact that
$\mu(S)^2 - 2 \int_{S} \po(x)d\mu(x)=\int_{S} (\mu(S)-2\po(x)) d\mu(x). $
\end{proof}

Let's now consider the function $c$ defined on $\Sigma$ for any $A$ of finite measure by 
$$c(A)= \int_{A} (\mu(A)-2\po(x)) d\mu(x).$$ 
The integral formula defined in \cref{iFo} means that for any selection $S$, we have $c(S)=\cho(S)$.
The second step of our proof of Theorem~\ref{theorem:po-solution} 
is to show that $T_k$ is the largest measurable set that maximizes 
$c$ for its respective measure.
More precisely, we will prove the following lemma.
\begin{myfle} \label{lemma-max-choice}
For any $k$ such as $\mu(T_k)< +\infty$ we have 
\begin{equation*}
c(T_k) 
	=\max\limits_{\substack{ A \in \Sigma, \mu(A) \leq \mu(T_k)}} c(A).
\end{equation*}
Moreover, if $\mu(A) \leq \mu(T_k)$ and $c(A) = c(T_k)$, then $A \subseteq T_k$ almost-everywhere.
\end{myfle}

The context of this
lemma is very similar to the knapsack problem \cite{martello1990knapsack}. 
In this problem, one needs to find a subset $A$ of a finite set of items $\{x_1, ...,x_n \}$ maximizing the total value $\Sigma_{x_i \in A} v(x_i)$ under the constraint that the total weight   
$\Sigma_{x_i \in A} w(x_i)$ of $A$ 
does not exceed a predefined maximum weight $w^*$. 
In Lemma~\ref{lemma-max-choice}, the total value is $c(A)$, the ratio of an element's value to its weight becomes $\mu(A)-2 \po(x)$, and the weight constraint is expressed as $\mu(A)\leq \mu(T_k)$.  
The solutions given by the lemma correspond to 
those yielded for the knapsack problem by George Dantzig's greedy approximation algorithm \cite{dantzig1957discrete}.
This algorithm consists of ordering elements by decreasing
value-to-weight ratio and then taking the $N$ first elements. 
$N$ is chosen in such a way, that taking one more element would cause excessive weight. 
The process of proving Lemma~\ref{lemma-max-choice} is similar to proving that George Dantzig's solutions are optimal for their respective weights.

\begin{proof}[Proof of Lemma~\ref{lemma-max-choice}] 
Let's consider $T_k$ such that $\mu(T_k)<+\infty$. To prove $c(T_k)$ is the maximum, we need to show for any $A\in \Sigma$, 
that $c(A) \leq c(T_k)$ if $\mu(A) \leq \mu(T_k)$.
As 
$A=(A \cap T_k) \cup (A \setminus T_k)$ and $T_k=(A \cap T_k)  \cup (T_k \setminus A)$,  requiring  $c(A) \leq c(T_k)$  is equivalent to requiring 
\begin{equation} \label{danzigineq}
\int \limits_{A \setminus T_k}(\mu(A)-2\po(x))d\mu(x) \leq \int \limits_{T_k \setminus A}(\mu(T_k)-2\po(x))d\mu(x).
\end{equation}
By definition of $T_k$, 
we have that 
$\po(x)> k$ for $x \in A \setminus T_k$, while $\po(x) \leq k$ for $x \in T_k \setminus A$.
Therefore,
\begin{align} \label{lemma:ineq1}
\int \limits_{A \setminus T_k}(\mu(A)-2\po(x))d\mu(x) &\leq \mu(A \setminus T_k)(\mu(A)-2k),\\
\shortintertext{and} \label{lemma:ineq2}
\int \limits_{T_k \setminus A}(\mu(T_k)-2\po(x))d\mu(x) &\geq \mu(T_k \setminus A)(\mu(T_k)-2k).
\end{align}
The constraint $\mu(A) \leq \mu(T_k)$  means that $\mu(A \setminus T_k) \leq \mu(T_k \setminus A)$.
Thus, 
we have 
\begin{equation}  \label{lemma:ineq3}
\mu(A \setminus T_k)(\mu(A)-2k) \leq \mu(T_k \setminus A)(\mu(T_k)-2k). 
\end{equation}
Finally, combining \cref{lemma:ineq3} 
with 
\cref{lemma:ineq1} and \cref{lemma:ineq2} guarantees the inequality in \cref{danzigineq} 
under the constraint $\mu(A) \leq \mu(T_k)$.

Let us now consider $A \in \Sigma$ such that $\mu(A) \leq \mu(T_k)$ and 
$c(A) = c(T_k)$. 
We now proceed to show that $A \subseteq T_k$ almost everywhere.
If $c(A)=c(T_k)$, the inequality in \cref{danzigineq} must be an equality. 
Under the assumption $\mu(A) \leq \mu(T_k)$, we again have the inequalities in
\cref{lemma:ineq1}, \cref{lemma:ineq2} and \cref{lemma:ineq3}, which must be equalities if
\cref{danzigineq} is an equality. 
However, because $\po(x)> k\text{ for }x \in A \setminus T_k$, the inequality in \cref{lemma:ineq1} can only become an equality if $\mu(A \setminus T_k)=0$.
\end{proof}

Now everything is in place to prove the theorem.
\begin{proof}[Proof of Theorem \ref{theorem:po-solution}]
Let us first prove 
that the at least $k$-Pareto optimal elements 
$T_k$ offer maximum choice. 
Lemma~\ref{lemma-max-choice} says that 
on $\Sigma$ the set 
$T_k$ maximises $c$ for its respective measure.
We assumed $T_k$ is a selection.
At the same time, for any selection $c=\cho$. 
It means that $T_k$ offers maximum choice for its respective measure.
Moreover, from Lemma~\ref{lemma-max-choice} and from $\mathcal{S} \subseteq \Sigma$ directly results that if a selection $A$ offers maximum choice for $\mu(T_k)$, then $A \subseteq T_k$ almost everywhere.
\end{proof}

The theorem does not guarantee uniqueness. 
For example, in the case of the relation $\leq$ on $\mathbb{R}$ and the Lebesgue measure, selections are all left-unbounded intervals and the choice of any selection is $0$.  
However, if for any selection
\begin{equation} \label{theo:uniq}
A \subset T_k \  \& \  \mu(A)<\mu(T_k) \implies 
\cho(A) < \cho(T_k),
\end{equation} 
then $T_k$ is a unique maximum. This is a direct
consequence
from the fact that $T_k$
contains
any other selection offering maximum choice for $\mu(T_k)<+\infty$.

Also, the proof of the theorem indicates that transitivity is not necessary. 
The set $T_k$ is only required to be a selection. 
This is the case for Lebesgue area measure and the non-transitive relation $R_{\llcorner '}$ defined on the unit square by $(x_1, x_2) R_{\llcorner '} (y_1, y_2)$ iff $\frac{y_1}{2} \leq x_1 \leq y_1$ and  $\frac{y_2}{2} \leq x_2 \leq y_2$.
The sets $T_k$ 
fulfill the requirement of being selections if $R^*$ satisfies the weakened transitivity condition 
$xR^*y \implies \po(x) \leq \po(y)$.

\begin{figure}
	\centering{
		\begin{tikzpicture}
\node [label=right:{$(1,1)$}] (e1t1) at (1.5,1) {$\bullet$}; 
\node [label=right:{$(1,1)$}] (e1t2) at (3.5,1) {$\bullet$};
\node [label=right:{$(1,0)$}] (e1b) at (2.5,0) {$\bullet$};
\node [label=left:{$(4,0)$}] (ea) at (0,0) {$\bullet$};
\draw[directed] (e1t2) -- (e1b);
\draw[directed] (e1t1) -- (e1b);

\coordinate (aea) at (0, 1);
\coordinate (e1t1b) at ($(e1t1)+(0,-1)$);	
\coordinate (e1t2b) at ($(e1t2)+(0,-1)$);	
\coordinate (down) at (0, -0.3);
\coordinate (up) at (0, 0.5);
\coordinate (s1) at ($0.5*(ea)+0.5*(e1t1b)+(down)$);
\coordinate (s2) at ($0.4*(ea)+0.6*(e1b)+(aea)+(up)$);
\coordinate (s3) at ($-0.2*(ea)+1.35*(e1t2b)+(aea)+(up)$);
\coordinate (s4) at ($-0.4*(ea)+1.55*(e1t2b)+(down)$);
\draw [rounded corners] (s1)-- (s2)-- (s3)--(s4);
\node [label=right:$S$] at ($(s4)-0.5*(down)$) {};
\end{tikzpicture}
		}
    \caption{$S$ is a selection that offers maximum choice but cannot be constructed from $T_k$. Point are encoded with their measure and $k$-Pareto optimality: $(\mu, \po)$.}
    \label{fig:non-po-selection-with-max-choice}
\end{figure}
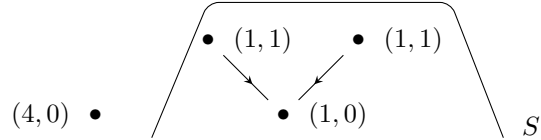

It is possible to prove that for
the counting measure, finite sets, and partial order relations,
all selections offering maximum choice can be obtained when sorting by
$k$-Pareto optimality and taking the first $l$ elements. 
However, for discrete measures with non-constant weights the above construction might not yield all selections offering maximum choice. 
A relative example is given in \cref{fig:non-po-selection-with-max-choice}. The partial order relation represented by its Hasse diagram. \cite{davey2002introduction}. 
The set $S$ is not of the form $T_k^* \cup A$ with $A \subseteq T_k\setminus T_k^*$. Nevertheless, it offers maximum choice for $m=3$.

\section{Further Theoretical Explorations}
\label{sec:theoretical-exploration}

\subsection{Efficient Computation of Choice}
\label{sec:theoretical-exploration:computation-of-choice}

Computation of choice may be simplified by performing the change of variable $y=\po(x)$ in the integral formula (\ref{iFo}). 
\begin{mypro}
\label{coro:simpl-choice-computation}
For any selection of at least $k$-Pareto optimal elements $T_k$ such that $\mu(T_k)<+\infty$
\begin{equation}\label{intfrorm:rk}
\cho(T_k) 
          =\mu(T_k)^2-2\int_{[0, k]} x d(\po * \mu)(x),
\end{equation}
where $\po * \mu $ is the image measure defined by $(\po * \mu)([a, b]) =\mu(\po^{-1}([a, b]))$. 
\end{mypro}

\subsection{Characteristics of Random Vectors}
\label{sec:theoretical-exploration:random-vectors}

Let us study a probability space $(\Omega, \Sigma, P)$.
We consider two random variables $X_1$ and $X_2$, as well as
the partial order relation $R_{\Omega \llcorner}$ defined on $\Omega$ as follows:
$$
\omega_x R_{\Omega \llcorner} \omega_y
   \text{ iff } 
X_1(\omega_x) \leq X_1(\omega_y)\text{ and }  X_2(\omega_x) \leq X_2(\omega_y)
.$$ 
For convenience, we use the following notations interchangeably:
$(X_1(\omega_x), X_2(\omega_x))$, $(x_1, x_2)$ or simply $x$.

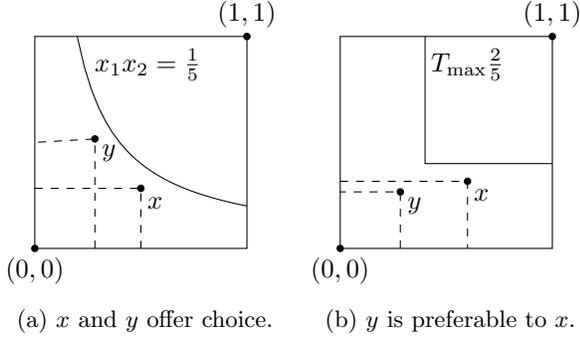
\begin{figure}
     \centering
     \begin{subfigure}[b]{0.23\textwidth}
         \centering
\begin{tikzpicture}[scale=.94]  
\coordinate (x) at(1.5,0.85);
\coordinate (y) at(0.85,1.55);
\coordinate (rect11) at (3,3); 
\coordinate (rect00) at (0,0); 
\coordinate (leftDown) at(0.3,-0.2);
\draw (rect00) to (3,0) to (rect11) to (0,3) to (0,0);
\draw [dashed](0,0.85) to (x) to (1.5,0);
\draw [dashed](0.85,0) to (y) to (0,1.5);
\fill (x) circle [radius=1.5pt];
\fill (y) circle [radius=1.5pt];
\fill (rect00) circle [radius=1.5pt];
\fill (rect11) circle [radius=1.5pt];
\draw [domain=0.6:3] plot (\x, {1.8/ \x} );
\node at ($(x)+(leftDown)+(-0.1, 0)$){$x$};
\node at ($(y)+(leftDown)+(-0.1, 0)$){$y$};
\node at ($(1.3, 2.8)+(leftDown)$){$x_1 x_2=\frac{1}{5}$};
\node [below] at (rect00){$(0,0)$};
\node [above] at (rect11){$(1,1)$};
\end{tikzpicture}
         \caption{$x$ and $y$ offer choice.}
         \label{fig:hyperbola:choice}
     \end{subfigure}
     \begin{subfigure}[b]{0.23\textwidth}
         \centering
\begin{tikzpicture}[scale=.94]  
\coordinate (x) at(1.8,0.95);
\coordinate (y) at(0.85, 0.8);
\coordinate (rect11) at (3,3); 
\coordinate (rect00) at (0,0); 
\coordinate (leftDown) at(0.3,-0.2);
\draw (rect00) to (3,0) to (rect11) to (0,3) to (0,0);
\draw [dashed](0,0.95) to (x) to (1.8,0);
\draw [dashed](0.85,0) to (y) to (0,0.8);
\fill (x) circle [radius=1.5pt];
\fill (y) circle [radius=1.5pt];
\fill (rect00) circle [radius=1.5pt];
\fill (rect11) circle [radius=1.5pt];
\draw  (1.2, 3) to (1.2, 1.2) to (3, 1.2);
\node at ($(x)+(leftDown)+(-0.1, 0)$){$x$};
\node at ($(y)+(leftDown)+(-0.1, 0)$){$y$};
\node at ($(1.5, 2.8)+(leftDown)$){$T_{\max} \frac{2}{5}$};
\node [below] at (rect00){$(0,0)$};
\node [above] at (rect11){$(1,1)$};
\end{tikzpicture}
         \caption{
         $y$ is preferable to $x$.
         }
         \label{fig:hyperbola:no-choice}
     \end{subfigure}
        \caption{Characterization of hyperbola. 
        }
        \label{fig:hyperbola}
\end{figure}

In this case, the selections are the sets situated below any
decreasing curve,
like for example the hyperbola $x_1 x_2 = 1 / 5$ and the curve defined by $max(x_1, x_2) = 2/ 5$ in Fig.~\ref{fig:hyperbola}.
The $k$-Pareto optimality $\po((x_1, x_2))$ is the joint cumulative probability distribution function $F(x_1, x_2)$.
Finally, the at least $k$-Pareto optimal elements are situated below the curve $F(x_1, x_2)=k$. 
Having choice between $x$ with coordinates $(x_1, x_2)$ and $y$ with coordinates $(y_1, y_2)$ means the 
rectangles $((0, 0),(x_1, 0),x,(0, x_2))$ and  $((0, 0),(y_1, 0),y,(0, y_2))$ are not 
nested, as depicted in Fig.~\ref{fig:hyperbola:choice}.

\subsubsection{Continuous Independent Variables}
\label{sec:theoretical-exploration:random-vectors:joint-prob}

If $X_1$ and $X_2$ are continuous, then the function $\po$ can be considered to be a linear extension of $R_{\Omega  \llcorner}$.
Indeed, if $x$ and $y$ are
chosen independently and at random from $\Omega$, then $P(\po(x)=\po(y))=0$.
On the other hand, it can be shown that  $\divers(\Omega)=1/2$.
This means that in half of the cases the relation $R_{\Omega \llcorner}$ 
cannot tell which element is preferable out of two elements extracted independently and at random from $\Omega$. 

Moreover if $X_1$ and $X_2$ are independent, then the condition for uniqueness from \cref{theo:uniq} holds. 
In the special case where $X_1$ and $X_2$ are uniformly distributed on the interval $[0,1]$, the joint cumulative probability distribution function takes form of $F(x_1, x_2) = x_1 x_2$ and the at least $k$-Pareto optimal elements are situated below a hyperbola $x_1 x_2 = k$.
Combining this fact with  Corollary~\ref{coro:max-div}  yields a surprising characterization of hyperbola, see \cref{fig:hyperbola} and the statement below.

\begin{mychar}
Out of all descending functions $f$ from $[0,1]$ to $[0,1]$ delimiting an area $\int_0^1 f(x) dx =c$, the arc of hyperbola is the one offering the highest likelihood the rectangles $((0, 0),(x_1, 0),x,(0, x_2))$ and  $((0, 0),(y_1, 0),y,(0, y_2))$ are not nested for two points $x$
and $y$
being drawn independently and at random from the delimited area.
\end{mychar}

\subsubsection{Independence on Marginal Distribution}
\label{sec:theoretical-exploration:random-vectors:marginal-independ}

Sklar's theorem \cite{sklar1959fonctions, durante2013topological} states that the cumulative distribution function $F(x_1, x_2)$ can be represented as $C(F_1(x_1), F_2(x_2))$ for a copula $C$.
Marginas of $X_1$ and $X_2$ are fully described by the marginal cumulative probability distributions $F_1$ and $F_2$, whereas the copula describes the dependence structure between $X_1$ and $X_2$.
The copula can be considered as a joint cumulative distribution function having two uniform marginal distributions on $[0, 1]$. 
Below we show that the introduced concepts 
do not depend on the marginal distribution of $X_1$ and $X_2$.

\begin{mypro}
For a continuous random vector $(X_1, X_2)$ and the relation
$R_{\Omega \llcorner}$, $k$-Pareto optimality, choice, diversity and selections offering maximum choice only depend on the copula $C$ of $X_1, X_2$.
\end{mypro}

\begin{proof}
Let us consider the mapping 
$$\begin{array}{ccccc}
G & : & \Omega & \to & [0,1]^2, \\
 & & \omega_x & \mapsto & (F_1(x_1), F_2(x_2)).
\end{array}$$
We consider $R_{\Omega \llcorner}$ and $P$ defined on $\Omega$. 
At the same time, on $[0,1]^2$ we consider $R_{\llcorner}$ defined by
$(x_1, x_2) R_{\llcorner} (y_1, y_2)$ iff $x_1 \leq x_2$ and $y_1 \leq y_2$, as well as the image measure $G*P$ defined 
on $[0,1]^2$ by $(G*P) (A)= P(G^{-1})(A)$. 
The map $G$ preserves probabilities in the sense that for any measurable 
$A$ in $\Omega$ we have $(G*P)(G(A))=P(A)$.
Moreover, $G$ preserves the relations in the sense that 
$x R_{\Omega \llcorner} y$ iff $G(x) R_{\llcorner} G(y)$. 
Selections are preserved in the sense that if $S$ is a selection for $R_{\Omega \llcorner}$, then $G(S)$ is a selection for $R_{ \llcorner}$. 
Ignoring negligible subsets, this mapping between selections is one-to-one. 
Therefore, $G$ also preserves selections, $k$-Pareto optimality, choice, and diversity.

The proposition finally results from the fact that  $G*P$  only depends on the copula. This is a consequence of
the fact that for any $(a_1, a_2)\in [0,1]^2$, we have 
$(G * P) ([0, a_1] \times[0, a_2])= C(a_1, a_2)$. This equality is a result of the following statements: 1) continuity which guarantees that $a_1$ and $a_2$ can be written as $F_1(x_1)$ and 
$F_2(x_2)$ for some appropriate $x_1$ and $x_2$; 2) 
the definition of image measure;
3) the fact that $G=(F_1, F_2) \circ (X_1, X_2)$; and 4) the equality $F(x_1, x_2)=C(F_1(x_1), F_2(x_2))$.
\end{proof}


\begin{mypro}
\label{pro:margin-dist}
If $X_1$ and $X_2$ are two continuous independent random variables, then for the relation 
$R_{\Omega \llcorner}$
\begin{align*}
P(T_k)&=k-k\ln(k),\\
\cho(T_k)&=(k-k\ln(k))^2-k^2\left( \frac{1}{2} - \ln k \right).
\end{align*}
\end{mypro}

\begin{proof}
Let us again consider the map $G$.
The image measure $G*P$ induced by $G$ on $[0,1]^2$ is the Lesbegue area measure. 
Independently of $P$, 
we have 
$G(T_k)=\{(x_1, x_2) \in [0,1]^2| x_1 x_2 \leq k\}$. 
Simple integration for \cref{intfrorm:rk} yields
\begin{align*}
P(\po^{-1}&(]-\infty, x])) =P(T_x), \\
                          &=\int_{(x_1, x_2) \in [0,1]^2}(G(T_x)) = x-x\ln(x). 
\end{align*} 
Therefore, $\po * P=-\ln(x) dx$ 
and \cref{intfrorm:rk} results in 
\begin{align*}
\cho(T_k)
         &=(k-k\ln(k))^2-2\int_0^k x(-ln(x))dx, \\
         &=(k-k\ln(k))^2-k^2\left( 1/2 - \ln k \right).    
\end{align*}
\end{proof}

\begin{figure}
	\centering{
	    \setlength{\figH}{7.5cm}
		\input{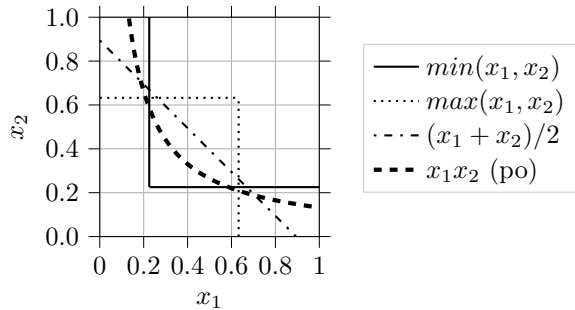}
		}
    \caption{
    Selections of the best 40\% according to different sorting criteria.
    }
    \label{fig:diversity_vs_area} 
\end{figure}

Now, it is possible to show that
$\lim_{k \to 0} \divers(T_k)=1.$
The fact that diversity slowly tends to the maximum possible value of 1 as $k$ tends to zero becomes even more surprising if one looks at the selections
\begin{align*}
T_{\min a}&=\{(x_1, x_2) \in [0,1]^2| \min(x_1, x_2) \leq a\}, \\
T_{\max a}&=\{(x_1, x_2) \in [0,1]^2| \max(x_1, x_2) \leq a\}.
\end{align*}
Here, 
$\lim_{a \to 0} \divers (T_{\min a})=\frac{3}{4}$, and
$\divers (T_{\max a})=\frac{1}{2}$. 
To further illustrate this observation, let us consider selections of a fixed measure $m$, which represent the best $m * 100$\% of elements, defined with different sorting criteria: minimum  $min(x_1, x_2)$, maximum  $max(x_1, x_2)$, average  $\frac{x_1 + x_2}{2}$, 
and Pareto optimality $x_1 x_2$\footnote{
Here, $\po$ is functionally related to geometric mean.
}.
As we can see in \cref{fig:diversity_vs_area},
$\min(x_1,x_2)$ delimits selections containing too many large values,
that is, extremes are overvalued.
On the other hand, 
all other sorting criteria except $\po$
undervalue extremes and include too many elements situated around the diagonal $x_1=x_2$.

Finally, in \cref{fig:diversity_vs_area} we demonstrate how diversity of the selections defined above depends on the fraction of selected elements $m$\footnote{
Note that in the considered case
for any selection $S$,  $\cho(S)=dx_1 dx_2(S)^2-\int_S\po(x)dx_1 dx_2$.
}. 
We can see that diversity is the largest when sorting by $k$-Pareto optimality.
This is a direct consequence of Theorem~\ref{theorem:po-solution}.

\begin{figure}
	\centering{
		\input{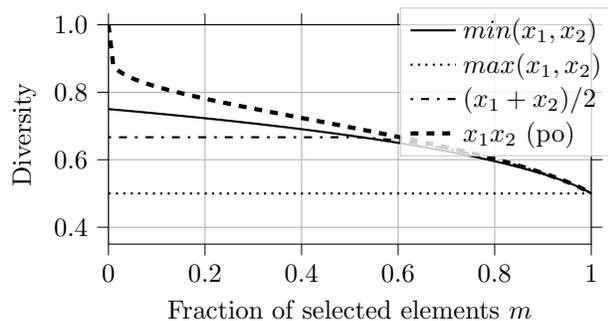}
		}
    \caption{
    Fraction of the best elements and diversity.
    }
    \label{fig:diversity_vs_area} 
\end{figure}

\section{Further Practical Explorations}
\label{sec:practical-exploration}

\subsection{Computation Complexity}
\label{sec:practical-exploration:computation-complexity}

A common solution for 
ranking $n$ elements of a set $X$ according to a partial order relation $R$
is to rank the elements according to their average ranking with respect to all linear extensions of $R$.
However, the total number of linear extensions exponentially increases with $n$, and the resulting algorithms are complex and slow \cite[p. 48]{loof2010gentphd}. 
For example, random sampling of linear extensions has an expected running time of $O(n^3 \log n)$ \cite{ huber2006sampling}. 
Below we show that sorting by $k$-Pareto optimality offers an efficient alternative.


The basic algorithm for the $k$-Pareto optimality based sorting is straightforward.
In the case of an arbitrary relation $R$, $\po(x)$ is computed by summing up the measures of the items that are strictly preferable to $x$.
This requires one pass through the whole set $X$ for every element $x \in X$ with computation complexity $O(n^2)$. 
The complexity of sorting $X$ by increasing values of $\po$ is $O(n \log n)$.
Therefore, the total computational complexity is $O(n^2)$.   



The case of composite relations defined on the probability space allows constructing even faster sorting procedures.
We illustrate this idea for $R=R_l\&R_p\&R_r$ from our housings example. 
We define the component relations as follows: for $i\in \{l,r,p\}$, $a R_i b$ iff 
$X_i(a) \leq X_i(b)$, where the real valued random variable $X_l$ represents proximity to the location, $X_p$ represents population size, and $X_r(x)$ is $0$ when $x$ is close to a river and $1$ otherwise.
For independent $X_i$, we have: 
\begin{align*}
&\po(x)=P(\{y|yR^*x\}),   \nonumber \\
    &=P(\{y|y R x\}) - P(\{y|y R x\myand xRy\}), 
                                                             \nonumber \\
      &=\prod_{i \in \{1,\dots,m\}}P( X_i \leq x_i) 
           - \prod_{i \in \{1,\dots,m\}}P( X_i = x_i). 
\end{align*}
The cumulative probability distributions $F_i(x)=P( X_i \leq x)$ can be approximated by the respective empirical cumulative probability distributions $\hat{F}_i(x)$.
The computation complexity of estimating $\hat{F}_i$ 
is $O(n \log n)$.
This needs to be done for every component relation $R_i$, resulting in the total complexity of $O( n \log n)$.


\subsection{Application to Genetic Optimization}
\label{sec:practical-exploration:genetic-algos}

In \cref{sec:introduction}, we hypothesized that sorting with choice maximization can be beneficial for genetic optimization.
Indeed, this strategy results in the maximization of the population diversity and allows exploring the search space more efficiently.
Additionally, Pareto dominance-based many-objective\footnote{
Concerns problems with 4 and more objectives.
} genetic optimization algorithms are known to
suffer from the lack of selection pressure \cite{palakonda2018pareto}. 
When the number of objectives increases, the number of
incomparable solutions grows exponentially.
However, as shown in \cref{sec:theoretical-exploration:random-vectors}, sorting random independent vectors by their Pareto optimality can be considered as a liner extension of the defined preference relation.
The fact that $P(\po(x)=\po(y))=0$ means that such sorting rarely produces ties and for any two solutions either $x$ is preferable to $y$ or vice versa.
In the rest of this subsection, we demonstrate that the proposed approach indeed improves the performance of genetic algorithms in the case of independent objectives.

To evaluate the proposed sorting procedure, we use it in \textit{NSGA-II} instead of
Pareto dominance-based sorting.
We experiment with two measures $\mu$: counting and probability measures.
This gives us two versions of genetic algorithms 
referred to as \textit{PO-count} and \textit{PO-prob} respectively.
These algorithms are compared with
implementations of the state-of-the-art algorithms 
\textit{NSGA-II} and \textit{NSGA-III}~\cite{deb2013evolutionary}
from the \textit{deap} python library\footnote{\url{https://deap.readthedocs.io/en/master/}}.
For the experimental evaluation, we use the $0/1$ knapsack problem
with independent objectives as defined  in~\cite{zitzler1999multiobjective}.
The number of knapsacks (objectives) is varied within the following set
$n_k \in \{2-8,10,15,25\}$ and the number of items is set to 250.
We adopt random selection with replacement and uniform crossover
with mutation probability 0.01.
We set the population size to 250 and the number of generations to 500.
All results are the average among 30 independent runs.

\begin{figure}
    \centering
\begin{tikzpicture}

\begin{axis}[
legend cell align={left},
legend style={fill opacity=0.8, draw opacity=1, text opacity=1, at={(0.03,0.97)}, anchor=north west, draw=white!80!black},
tick align=outside,
tick pos=left,
x grid style={white!69.0196078431373!black},
xlabel={Number of knapsacks, \(\displaystyle n_k\)},
xmajorgrids,
xmin=0.85, xmax=26.15,
xtick={1,2,3,4,5,6,7,8,10,15,25},
xtick style={color=black},
y grid style={white!69.0196078431373!black},
ylabel={Hypervol. increase,\%},
ymajorgrids,
ymin=-15.8234209671474, ymax=64.9014174796205,
ytick={-10, 0,10,20,30,40,50, 60},
ytick style={color=black},
height=4.5cm,
width=8.cm,
]
\addplot [very thick, black]
table {%
2 -0.400450982727018
3 -0.73294126551055
4 -1.94921885533832
5 -7.10858286755829
6 -10.6682504724551
7 -12.154110128658
8 -11.4512546306928
10 -10.7822381254516
15 -4.23483441969206
25 -0.0369180048387656
};
\addlegendentry{NSGA-III}
\addplot [very thick, black, dashed]
table {%
2 -0.1912121653585
3 -0.804234096148512
4 -0.132760059848866
5 -0.145387105898142
6 0.331737466914415
7 -0.567785946624202
8 0.226926279581223
10 -0.314466908663837
15 -0.139531721123252
25 0.0732254456485776
};
\addlegendentry{PO-count}
\addplot [very thick, black, dotted]
table {%
2 4.35522898503701
3 3.17869364747422
4 2.15259143928172
5 -1.35546644474396
6 -1.62902327235098
7 -1.13987592611693
8 0.760690187050005
10 6.06715884796208
15 24.1341824073453
25 61.232106641131
};
\addlegendentry{PO-prob}
\end{axis}

\end{tikzpicture}
    \caption{Increase in hypervolume compared to \textit{NSGA-II}.}
    \label{fig:hypervolume-increase}
\end{figure}
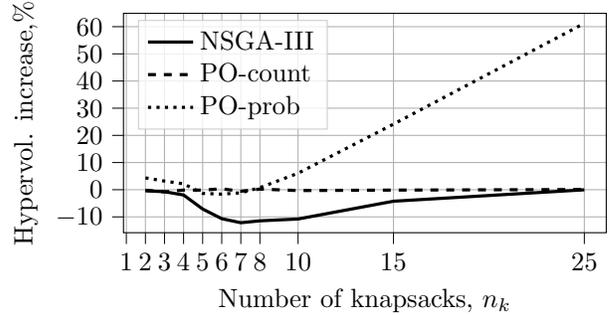

Below we analyze the performance of different algorithms in terms of the classical \textit{hypervolume} metric \cite{shang2020survey} with the origin of coordinates as a reference point.
In our setup, this metric is to be maximized.
We choose \textit{NSGA-II} as the baseline,
and present the relative changes in the \textbf{hypervolume} indicator for the rest of the algorithms  in \cref{fig:hypervolume-increase} (increase: positive number, decrease: negative number).
We notice that despite having been developed for the many-objective optimization, \textit{NSGA-III} almost always results in lower values of hypervolume, even for a large number of knapsacks.
This confirms a similar observation from~\cite{ishibuchi2016performancega}, 
and supports our choice of \textit{NSGA-II} as a baseline for implementation and comparison instead of \textit{NSGA-III}.
Further, we see that the value of relative increase for \textit{PO-count} is always very close to $0$.
It means that \textit{PO-count} yields a population covering the same hypervolume as \textit{NSGA-II}.
Contrarily, \textit{PO-prob} improves the hypervolume, as compared to \textit{NSGA-II}.
This difference is visible for small $n_k$ (+4\% for $n_k=2$)
and is especially prominent
for large $n_k$ ($+60\%$ for $n_k=25$).
For $n_k$ between 5 and 7, \textit{PO-prob} results in lower values of hypervolume
than \textit{NSGA-II}.
However, the relative decrease
in thes cases 
does not exceed $-1.63$\%. 
Also, within this range, \textit{PO-count} performs slightly better than other algorithms.
These results demonstrate that the proposed approach improves the performance of genetic algorithms, especially in the case of many-objective optimization. 
It also suggests that the choice of the measure $\mu$ has a large impact on the performance.
The latter relationship will be studied in future work.

\section{Conclusion}
\label{sec:conclusions}

In this paper, we formulate the problem of \textit{generalized topological sorting with choice maximization}, which, to the best of our knowledge, was not considered in the literature before.
We also prove
that the \emph{at least $k$-Pareto optimal sets} provide unique solutions.
Further theoretical analysis of this problem leads us to an interesting relationship between the diversity of random points and the arc of hyperbola. 
Additionally, 
we propose a computationally efficient algorithm for the calculation of $k$-Pareto optimality for 
probability measures.
Finally, we demonstrate a successful application of the developed theory.
We show that sorting by $k$-Pareto optimality can drastically improve the performance of many-objective genetic optimization algorithms.
In our experiments, the proposed solution based on the probability measure allows increasing the value of hypervolume by up to 60\% for 25 objectives.
This result can be considered as a potential solution to the problem of searchability deterioration in Pareto-dominance optimization.

We also believe that the proposed general framework can be used in different applications. 
In future work, we plan to study the applicability of $k$-Pareto optimality for constrained optimization, scheduling problems, recommender systems, and the development of statistical indicators.
Maximization of choice might be also useful when studying causality and fairness.


\bibliography{literature}  

\onecolumn
\aistatstitle{
Supplementary Materials}

\appendix

\section{Further Examples for Simple Relations and Measures}
\label{annex:further-examples}
In the first two supplementary examples, we consider a situation typical in economics or multi-objective optimization. 
Later, we show how the proposed concepts apply to arbitrary 
transitive relations. 
  

\subsection{Continuous Measures}
\label{annex:further-examples:continuous-measures}

Let us again consider $R_{\llcorner}$ as defined in \cref{sec:pareto-optimality:examples}. The relation 
$R_{\llcorner}$ models preference for small values of $x_1$ and $x_2$. However, instead of
assuming $X$ to be a finite subset of $\mathbb{R}^2$, we now study the unit square $[0,1]\times[0,1]$ with three continuous measures: 
the Lebesgue area measure $d x_1 d x_2$, as well as $2x_2 d x_1 d x_2$ and $4 x_1 x_2 d x_1 d x_2$. 
In each case, the total measure of the unit square equals to one. 
The Lebesgue area measure represents 
elements with two uniformly distributed characteristics $x_1$ and $x_2$;
$2x_2 dx_1 dx_2$ represents rarefaction of items having small values of $x_2$, whereas $4 x_1 x_2 dx_1 dx_2$ represents rarefaction of items having small values of both $x_1$ and $x_2$.

\begin{figure}[!h]
     \centering
     \begin{subfigure}[b]{0.48\textwidth}
         \centering
         \input{plots/rarefication_dx1_dx2}
         \caption{
         At least $k$-Pareto optimal elements of measure 0.1 for the relation $R_{\llcorner}$ and three measures: the Lebesgue area measure $d x_1 d x_2$,  
         rarefaction of $x_2$ defined by $2 x_2 d x_1 d x_2$, and rarefaction of $x_1$ and $x_2$ defined by 
     $4 x_1 x_2 d x_1 d x_2$.
         }
         \label{subfig-1:rarification}
     \end{subfigure}
     \hfill
     \begin{subfigure}[b]{0.48\textwidth}
         \centering
         \input{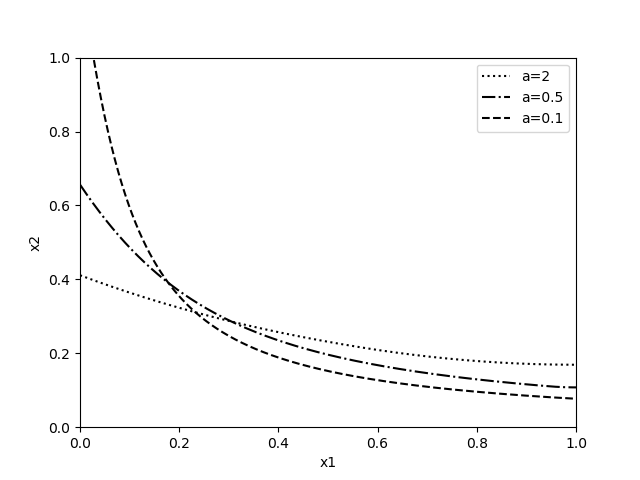}
         \caption{
         At least $k$-Pareto optimal elements of measure 0.25 for the cone-based relations $R_a$ and
         the Lebesgue area measure.
         Larger values of $a$ represent higher maximum accepted trade-offs.
         }
         \label{subfig-2:tradeoff}
     \end{subfigure}
        \caption{Examples at least $k$-Pareto optimal elements with two continuously distributed characteristics $x_1$, $x_2$.
        }
        \label{fig:dummy}
\end{figure}

For each of the cases defined above,
in \cref{subfig-1:rarification}
we show
the set of at least $k$-Pareto optimal elements of measure $0.1$, which corresponds to selecting the $10$ best percent.
All three sets demonstrate the qualitative behaviours 
expected from sets delimited by indifference curves when the corresponding rarefaction occurs. 
Indeed, the curve corresponding to the uniform distribution and the Lebesgue area measure $dx_1 dx_2$ is symmetric. Also, in this case,
$\po(x_1, x_1)= x_1 x_2$, and the sets of at least 
$k$-Pareto optimal elements are the sets situated below arcs of hyperbola defined by the equation $x_1 x_2=k$, see \cref{sec:theoretical-exploration:random-vectors,sec:theoretical-exploration:random-vectors:joint-prob} for more details.
Applying rarefaction with respect to $x_2$ prioritises smaller values of this characteristic.
This is represented by shifting upwards the right part of the hyperbola arc, see the curve for $2x_2 dx_1 dx_2$.
Indeed, in this case, the small values of $x_2$ are observed less often.
This results in selecting additional elements with large values of $x_1$ but relatively small values of $x_2$ to compensate for this rarefaction.
Finally, rarefaction with respect to both $x_1$ and $x_2$ results in the fact that the small values of both characteristics are observed less often.
Thus, elements with larger values of $x_1$ and $x_2$ should be selected to generate a selection of the required measure.
It results in the shift of the hyperbola upwards following the direction of the main diagonal, see the curve for $4 x_1 x_2 d x_1 d x_2$.

   
\subsection{Cone-based Relations}
\label{annex:further-examples:cone-relation}

\begin{figure}[!h]
     	\centering{
       	\tikzmath{\ssl = 3;  
 \x1 =1.5; \x2=\x1 +0.5; 
 \h = 1.4;
 \excess=0.25;
 \fillepsilon=0.1; 
 \tradeoff=1; 
 \tradeoffI=(\tradeoff / \x2)*(\ssl - \x2); 
 \hI=(\tradeoff / \x2)*\h; 
 \hExcess=(\excess / \x2)*\h;
 } 
 
\begin{tikzpicture}
\coordinate (x) at(\x1,\x2);
\coordinate (rect11) at (\ssl,\ssl); 
\coordinate (rect10) at (\ssl,0); 
\coordinate (rect01) at (0,\ssl); 
\coordinate (rect00) at (0,0); 
\coordinate (leftDown) at(0.3,-0.2);
\coordinate (tradeoffDown) at ($(\x1,0)+(\tradeoff,0)$); 
\coordinate (tradeoffUp) at ($(\x1,\ssl)-(\tradeoffI,0)$);
\coordinate (hCorner) at ($(x)-(0,\h)$);
\coordinate (hIntercept) at ($(hCorner)+(\hI, 0)$);
\fill[pattern=horizontal lines, pattern color=gray!50]       
($(0,\ssl)+(-\excess,-\fillepsilon)$) 
-- ($(0,\x2)+(-\excess,\fillepsilon)$) 
-- ($(x)+(-\fillepsilon, \fillepsilon)$)
--($(tradeoffUp)+(-\fillepsilon,-\fillepsilon)$) ;

\fill[pattern=horizontal lines, pattern color=gray!50]       
($(x)+(\fillepsilon,-\fillepsilon)$) 
-- ($(tradeoffDown)+(\fillepsilon,0)+(\hExcess,-\excess)$) 
-- ($(rect10)+(-\fillepsilon, \hExcess)$)
--($(\ssl,\x2)+(-\fillepsilon,-\fillepsilon)$) ;

\draw (rect00) to (rect10); 
\draw (rect00) to (rect01);
\draw  [->](0,0) coordinate (Oo) -- (0,\ssl) coordinate (Ox) ;
\draw  [->](Oo) -- (\ssl,0) coordinate (Oy) ;
\draw (-\excess,\x2) to (x) to ($(tradeoffDown)+(\hExcess,-\excess)$);
\draw [dashed](tradeoffUp) to (x) to (\ssl,\x2);

\draw [dotted] (x) to (hCorner) to (hIntercept);
\draw[decorate,decoration={brace, aspect=0.3, amplitude=4, raise=1.5}] (hCorner) --  (x) node[pos=0.3,label=left:{\footnotesize $h$}]{} ;
\draw[decorate,decoration={brace, aspect=0.5,  amplitude=4, raise=1.5}]  (hIntercept) 
--  (hCorner) node[pos=0.5,label=below:{\footnotesize $a h$}]{} ;
\fill (x) circle [radius=1.5pt];
\fill (rect00) circle [radius=1.5pt];
\node [label=above right:{$x$}] at (x) {};
\node [label=below right:{$y R^{*}_a x$}] at (0,\x2) {};
\node [label=below:{$x R^{*}_a y$}] at ($(rect11)+(-0.3,0)$) {};
\node [left] at (rect00){$(0,0)$};

\node [label=left:{$x_2$}] at (rect01) {};
\node [label=below:{$x_1$}] at (rect10) {};
\end{tikzpicture}
     }
     \caption{An illustration of a cone-based relation $R_a$.}
     \label{subfig-2:R_a}
  \end{figure}
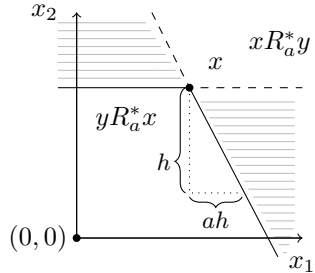

Let us again consider the unit square, the Lebesgue measure, and a positive constant $a$. 
However, this time the preference relation is defined as follows: $y R_a x$ iff $y_2 \leq x_2$ and $x_2-y_2 \geq a (y_1-x_1)$. The above relation $R_a$ is an example of a cone-based relation illustrated in \cref{subfig-2:R_a}.
This relation has the intuitive meaning of
giving up 
$(x_1, x_2)$ for getting $(y_1, y_2)$ if the improvement (diminution) in the second characteristic is at least $a$ times the trade-off (increase) in the first characteristic. In this case, selections are sets 
delimited by descending curves $x_2=f(x_1)$ such that
$-\frac{1}{a} \leq\frac{d f}{d x_1}\leq 0$.

Let us now consider the sets of at least $k$-Pareto optimal elements of measure $0.25$ for the three values of $a$: $a=\frac{1}{10}$, $a=\frac{1}{2}$, and $a=2$, see \cref{subfig-2:tradeoff}. 
Larger values of $a$ represent higher maximum accepted trade-offs.
This is represented by the gradual degeneration of the hyperbola into a straight horizontal line when $a$ increases.
As shown in the figure, the three sets demonstrate plausible behavior. 
In the situation discussed in \cref{annex:further-examples:continuous-measures},
the relation $R_{\llcorner}$ corresponds to the extreme case of the relation $R_a$ with $a=0$.

\subsection{Transitive relations}

\begin{figure}[!h]
	\centering{
		\begin{tikzpicture}
\node (i) at (0,0) [label=right:$6$] {$i$};
\node [below of=i](g2) {};
\node [left of=g2, label=right:$4$](g) {$g$};
\node [left of=g, label=right:$3$](f) {$f$};
\node [right of=g2, label=right:$2$](h) {$h$};
\node [below of=g, label=right:$2$](d) {$d$};
\node [below of=g2, label=right:$1$](e) {$e$};
\node [below of=d, label=right:$0$](a) {$a$};
\node [below of=e, label=right:$0$](b) {$b$};
\node [right of=b, label=right:$0$](c) {$c$};

\coordinate (p0) at($1.7*(a)-0.7*(b)$);
\coordinate (p1) at($0.35*(f)+0.65*(d)$);
\coordinate (p2) at($0.4*(g)+0.6*(d)$);
\coordinate (p3) at($0.6*(g)+0.4*(e)$);
\coordinate (p4) at($0.5*(i)+0.5*(h)+(0.3,0.3)$);
\coordinate (p5) at($0.5*(i)+0.5*(h)+(0.8,0.4)$);
\draw [rounded corners] (p0) -- (p1)-- (p2)-- (p3)--(p4)--(p5);
\node (p5l) at ([shift=({0.2, -0.2})]p5) {$T_2$};

\draw[directed] (i) -- (g);
\draw[directed] (i) -- (h);
\draw[directed] (f) -- (d);
\draw[directed] (g) -- (d);
\draw[directed] (g) -- (e);
\draw[directed] (h) -- (e);
\draw[directed] (e) -- (b);
\draw[directed] (d) -- (b);
\draw[directed] (d) -- (a);
\end{tikzpicture}
		}
    \caption{An illustration of simple partial order relation. The values of $\po$ are shown by numbers, and the set of at least $2$-Pareto optimal elements $T_2$ is delimited by a curve.}
    \label{fig:hasse-diagram-po}
\end{figure}
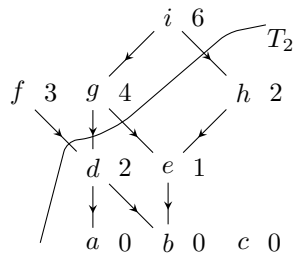

In general, it is possible to show that if $R^*$ is transitive, then for any $k$, the set $T_k$ is a selection.
In particular, if $R$ is a partial order relation,
$\mu$ is strictly positive, and $X$ is countable, we obtain a linear extension \cite{dantzig1957discrete} of $R$ when sorting $X$ by increasing values of $\po$ and sorting ties in any order. 
Selections are represented by downsets.
The latter are obtained when topologically sorting $X$ and taking the first $n$ elements, for any $n$. 
If $\mu$ is the counting measure, then $\po(x)$ is simply the number of elements that can be reached by following downwards the edges of the corresponding Hasse diagram. 
An example of such a relation represented by its Hasse  diagram is depicted in \cref{fig:hasse-diagram-po}.


\section{Further Theoretical Exploration of \em{at least k-Pareto optimal elements}}

\subsection{Sorting by $k$-Pareto Optimality versus Sorting by Pareto Fronts}

\begin{figure}
     \centering
     \begin{subfigure}[b]{0.48\textwidth}
         \centering
         \includegraphics[width=0.99\textwidth]{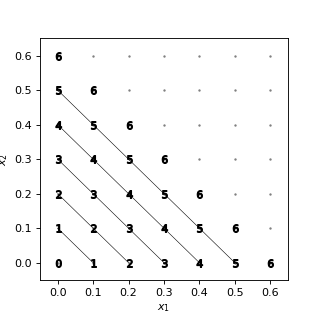}
         \caption{
         Sorting by Pareto fronts.
         }
         \label{fig:sorting-grid:front}
     \end{subfigure}
     \hfill
     \begin{subfigure}[b]{0.48\textwidth}
         \centering
         \includegraphics[width=0.99\textwidth]{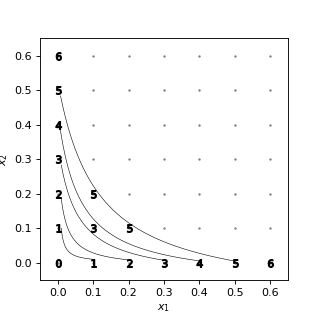}
         \caption{
         Sorting by $\po$.
         }
         \label{fig:sorting-grid:po}
     \end{subfigure}
        \caption{
        Sorting points of a grid. The equivalence classes (fronts) are represented by numbers.
        }
        \label{fig:sorting-grid}
\end{figure}

\begin{figure}[!h]
    \centering{
       \includegraphics[width=0.45\textwidth]{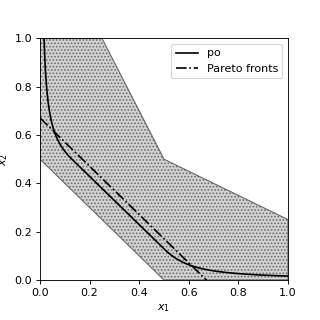}
       }
     \caption{
     Selections of $\mu = 0.2$ for the set $X$ composed of points in the shaded area.
     }
    \label{fig:PO_whiskers}
   \end{figure}

To further illustrate the difference between sorting by Pareto fronts and $k$-Pareto optimality, we demonstrate the results of sorting the elements of a 2-dimensional grid for the relation $R_{\llcorner}$ and counting measure in \cref{fig:sorting-grid}. 
The numbers on the plots represent the equivalence class (front) to which a point was assigned by the relative sorting procedure.
As we can see, sorting by Pareto fronts splits the points by straight lines, see \cref{fig:sorting-grid:front}.
At the same time, sorting by $\po$ results in splitting by hyperbola-like curves, see \cref{fig:sorting-grid:po}. 
We can also notice, that extreme solutions\footnote{
\textit{Extreme solution} here means that a solution is very good according to one criteria and is bad according to another.
} are valued more when sorting by $\po$.
Indeed, most of the non-extreme solutions are pushed to further equivalence classes, as compared to sorting by Pareto fronts.
This characteristic of $\po$-based sorting is also clearly visible in \cref{fig:PO_whiskers}.
Here we present selections of $\mu = 0.2$ for the same relation and the set $X$ composed of a large number of points placed on a regular grid within the shaded area.


\begin{figure}
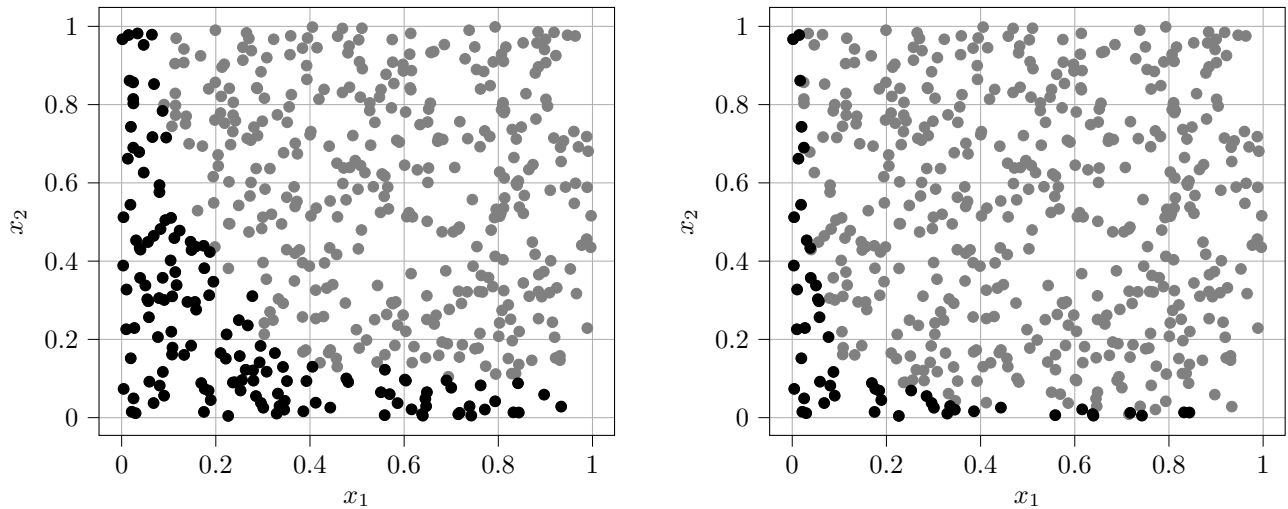

     \centering
     \begin{subfigure}[b]{0.48\textwidth}
         \centering
         \input{plots/sorting_uniform_front}
         \caption{
         Sorting by Pareto fronts. The total number of equivalence classes (fronts) is 39.
         }
         \label{fig:sorting-uniform:front}
     \end{subfigure}
     \hfill
     \begin{subfigure}[b]{0.48\textwidth}
         \centering
         \input{plots/sorting_uniform_po}
         \caption{
         Sorting by $\po$. The total number of equivalence classes (fronts) is 259.
         }
         \label{fig:sorting-uniform:po}
     \end{subfigure}
        \caption{
        Sorting 500 uniformly distributed points. Points in black belong to the first 10 equivalence classes (fronts). Note, that the total number of equivalence classes is larger for $\po$-based sorting. The latter approach results in fewer ties.
        }
        \label{fig:sorting-uniform}
\end{figure}

However, sorting by Pareto front does not always result in selections delimited by straight lines.
Analysing the results of sorting for uniformly distributed points, we observe that both sorting methods result in hyperbola-like selections, see \cref{fig:sorting-uniform}.
This means that sorting by Pareto fronts is more sensitive to the topological structure of the analyzed space, while sorting by $\po$ preserves its characteristics.

\subsection{Further Solutions of the Maximum Choice Problem}
\label{sec:theoretical-exploration:further-solutions-of-max-choice}



\subsubsection*{Further similar solutions}

It is possible to prove
that Theorem~\ref{theorem:po-solution} also holds for $T_k^*$ defined with a strict inequality ($<$) as follows, see Def.~\ref{def:po-set} for comparison.
    $$T_k^*=\{x \in X | \po(x) < k  \}.$$
In this case, the proof of the fact that $T_k^*$ is the largest in Lemma~\ref{lemma-max-choice} 
requires analysis of inequality~(\ref{lemma:ineq2}) instead of inequality~(\ref{lemma:ineq1}).
However, 
the proof of the fact that any selection $T$ such that 
$T_k^* \subseteq  T \subseteq T_k$  offers maximum choice for $\mu(T)$
becomes a bit more technical. Moreover,
it is possible to prove that $T$ is the largest selection of this kind.
Precisely, for any other selection $S$ offering maximum choice for 
$\mu(T)$, 
$S \subseteq T'$ for some $T'$ such that $\mu(T')=\mu(T)$ and $T_k^* \subseteq  T' \subseteq T_k$.

\subsubsection*{Completeness of the Solutions}
 
In the case condition in \cref{theo:uniq} holds, and if for any selection $S$ there is a $k$ and there are selections $T$ such that $T_k^* \subseteq  T \subseteq T_k$ and $\mu(T)=\mu(S)$,
then those selections $T$ are the only selections offering maximum choice.
Therefore, we have a complete list of selections offering maximum choice. 
This is the case for the typical example of the relation $R_{\llcorner}$ defined in \cref{sec:pareto-optimality:examples} and the Lebesgue area measure defined on the unit square $[0,1]^2$.
This also holds for any discrete measure with constant non-zero weights, 
for example, for the 
counting measure on a finite set with a partial order relation.
In the latter case, 
topologically sorting
by increasing values of $\po$ and then taking the first $n$ elements 
results in a set offering maximum 
choice.

\subsubsection*{Existence of Solutions of a Different Nature}

Let us consider again the example in \cref{fig:non-po-selection-with-max-choice}.
Here,
the set $S$ is not of the form $T_k^* \cup A$ with $A \subseteq T_k\setminus T_k^*$. Nevertheless, it offers maximum choice for $m=3$.




\subsection{Diversity}
\label{sec:theoretical-exploration:diversity}


As it was discussed in \cref{sec:pareto-optimality:max-choice}, the concept of \emph{maximum choice} is functionally related to the concept of \emph{diversity}, see Def.~\ref{def:diversity}.
In Theorem~\ref{theorem:po-solution} we cannot simply replace choice by diversity.
However, by considering only selections of a fixed measure, we obtain the following straightforward corollary of Theorem~\ref{theorem:po-solution}.

\begin{mycoro}
\label{coro:max-div}
If $R$ is a transitive relation,
then for any set of at least $k$-Pareto optimal elements 
$T_k$ such that $\mu(T_k)<+\infty$, we have
$$\divers(T_k) 
	=\max\limits_{\substack{ S \in \mathcal{S}, \mu(S) = \mu(T_k)}} \divers(S).$$ 
Moreover $T_k$ is the unique such maximum. Precisely if $S$ is a selection such that $\mu(S)=\mu(T_k)$, and $\divers(S)=\divers(T_k)$, then $S=T_k$ almost everywhere. 
\end{mycoro}

\section{Further Practical Exploration of \em{at least k-Pareto optimal elements}}

\subsection{Additional Results for Genetic Optimization}
\label{sec:practical-exploration:genetic-algos-further}

In \cref{sec:practical-exploration:genetic-algos}, we evaluate the performance of the genetic algorithms using the \textbf{hypervolume} indicator.
In this section, we further analyze the behavior of both the state-of-the-art and the proposed algorithms with respect to other metrics.
In particular, we study the fraction of solutions dominated by the solution of alternative algorithms and analyze the time complexity of the sorting procedure.

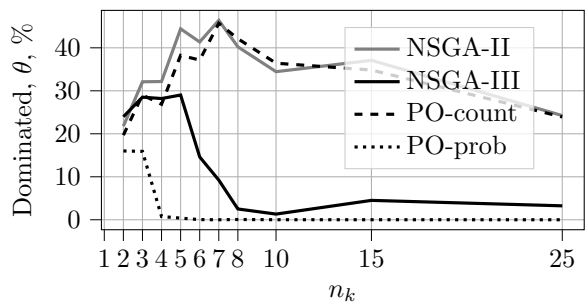
\begin{figure}[h!]
\centering
\begin{tikzpicture}

\begin{axis}[
legend cell align={left},
legend style={fill opacity=0.8, draw opacity=1, text opacity=1, at={(0.7,0.6)}, anchor=center, draw=white!80!black},
tick align=outside,
tick pos=left,
x grid style={white!69.0196078431373!black},
xlabel={\(\displaystyle n_k\)},
xmajorgrids,
xmin=0.85, xmax=26.15,
xtick style={color=black},
xtick={1,2,3,4,5,6,7,8,10,15,25},
y grid style={white!69.0196078431373!black},
ylabel={Dominated, $\theta$, \%},
ymajorgrids,
ymin=-2.3195, ymax=48.7095,
ytick style={color=black},
ytick={0,10,20,30,40},
height=4.5cm,
width=8.cm,
]
\addplot [very thick, white!50.1960784313725!black]
table {%
2 21.7702
3 32.0997
4 32.16
5 44.42
6 41.36
7 46.39
8 40.23
10 34.45
15 37.09
25 24.24
};
\addlegendentry{NSGA-II}
\addplot [very thick, black]
table {%
2 23.9707
3 28.5219
4 28.2
5 29.02
6 14.61
7 9.23
8 2.51
10 1.32
15 4.51
25 3.24
};
\addlegendentry{NSGA-III}
\addplot [very thick, black, dashed]
table {%
2 19.6855
3 29.0279
4 26.87
5 38.26
6 37.15
7 45.71
8 42.12
10 36.42
15 34.83
25 23.92
};
\addlegendentry{PO-count}
\addplot [very thick, black, dotted]
table {%
2 16.0072
3 15.909
4 0.7593
5 0.382
6 0.0306
7 0
8 0.06
10 0.01
15 0.01
25 0
};
\addlegendentry{PO-prob}
\end{axis}

\end{tikzpicture}
\caption{Average percentage of solutions dominated by other algorithms, $\theta$.}
\label{fig:domination}
\end{figure}

We calculate the percentage of \textbf{dominated solutions} as follows.
For a given pair of algorithms \textit{algorithm1} and \textit{algorithm2}, we calculate how many solutions of
\textit{algorithm2} (\textit{dominated algorithm}) are dominated by solutions of \textit{algorithm1}
(\textit{dominating algorithm}).
After that, we average the obtained results among all \textit{dominating algorithms} to get an average fraction
of dominated solutions, denoted by $\theta$.
Naturally, lower values of $\theta$ indicate better performance.
We present the corresponding results in \cref{fig:domination}.
We notice the following tendencies.
\textit{NSGA-II} and \textit{PO-count} behave very similarly.
For $n_k=2$, the value of $\theta$ for these algorithms is around 20\%.
After that, it starts increasing and reaches its peak of approximately 45\% for $n_k=7$.
Finally, it gradually decreases to 24\% for $n_k=25$.
\textit{NSGA-III} starts at a similar level and reaches its peak of approximately 30\% for $n_k=5$.
After that, it decreases below 10\% for $n_k=7$ and stays relatively close to 0 for the larger numbers of knapsacks.
These results demonstrate the superiority of \textit{NSGA-III} over \textit{NSGA-II}
in the case of many-objective optimization.
\textit{PO-prob} starts at around 16\%.
However, for $n_k=4$ the value of $\theta$ it already almost 0 and does not go up for larger numbers of knapsacks.
This shows that the solutions produced by this algorithm are rarely dominated.
Thereby, \textit{PO-prob} is an
effective approach for many-objective optimization problems.

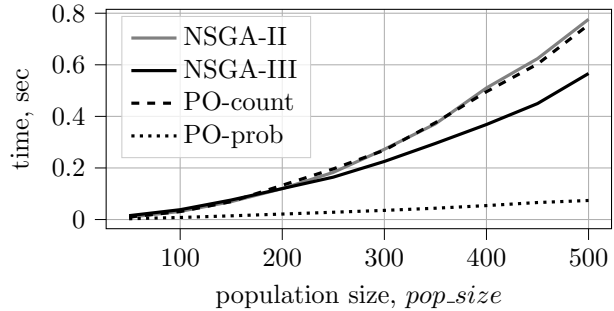
\begin{figure}
\centering
\begin{tikzpicture}

\begin{axis}[
legend cell align={left},
legend style={fill opacity=0.8, draw opacity=1, text opacity=1, at={(0.03,0.97)}, anchor=north west, draw=white!80!black},
tick align=outside,
tick pos=left,
x grid style={white!69.0196078431373!black},
xlabel={population size, $pop\_size$},
xmajorgrids,
xmin=27.5, xmax=522.5,
xtick style={color=black},
y grid style={white!69.0196078431373!black},
ylabel={time, sec},
ymajorgrids,
ymin=-0.0355372221469879, ymax=0.815024793863297,
ytick style={color=black},
height=4.5cm,
width=8.3cm,
]
\addplot [very thick, white!50.1960784313725!black]
table {%
50 0.00898474216461182
100 0.0315520668029785
150 0.0694082355499268
200 0.119809532165527
250 0.182883892059326
300 0.272381098270416
350 0.37105612039566
400 0.510234861373901
450 0.623849294185638
500 0.776362884044647
};
\addlegendentry{NSGA-II}
\addplot [very thick, black]
table {%
50 0.0149984502792358
100 0.0381470966339111
150 0.0756760239601135
200 0.120083999633789
250 0.164426438808441
300 0.225280463695526
350 0.295139102935791
400 0.368468396663666
450 0.449669146537781
500 0.566201233863831
};
\addlegendentry{NSGA-III}
\addplot [very thick, black, dashed]
table {%
50 0.00986402750015259
100 0.0308664488792419
150 0.0694935011863708
200 0.132789270877838
250 0.196425218582153
300 0.268778259754181
350 0.375608298778534
400 0.495700113773346
450 0.602707188129425
500 0.755122404098511
};
\addlegendentry{PO-count}
\addplot [very thick, black, dotted]
table {%
50 0.00312468767166138
100 0.00796788692474365
150 0.014371166229248
200 0.0212447810173035
250 0.0282719874382019
300 0.0354595541954041
350 0.0438962030410767
400 0.0540533208847046
450 0.0659357166290283
500 0.0738907241821289
};
\addlegendentry{PO-prob}
\end{axis}

\end{tikzpicture}
        \caption{Sorting duration as a function of population size for 10 knapsacks, $n_k=10$.}
        \label{fig:time:pop-size}
\end{figure}

In \cref{fig:time:pop-size}, we demonstrate the dependence of sorting time on the population size for values of $pop\_size$
ranging from 50 to 500.
The reported values are the averages over 100 independent executions of one iteration of the corresponding genetic algorithm.
From the figure, we can see that \textit{PO-prob} requires much less time than all other algorithm.
The results for \textit{NSGA-II} and \textit{PO-count} tend to be very close, as in other experiments.
This observation also has theoretical explanation.
Indeed, choosing the next generation for \textit{NSGA-II} and \textit{NSGA-III} has time complexity of $O(N^2M)$ and $max\{ O(N^2M), O(N^2 log_{M-2}N) \}$ 
respectively where $M$ stands for number of objectives and $N$ is the population size, see~\cite{deb2002fast,deb2013evolutionary}.
At the same time, sorting in \textit{PO-prob} comes down to independent sorting procedures with respect to
every objective.
The time complexity of this procedure is $O(NMlog(N))$.
These results are in line with the theoretical analysis presented in \cref{sec:practical-exploration:computation-complexity} and prove the computational efficiency of the approximate ranking calculation procedure used in \textit{PO-prob}.

The maximum choice theorem (Theorem ~\ref{theorem:po-solution}) has an intuitive interpretation in the context of genetic algorithms.
Assume that the selection step is required to pick a selection of a given maximum size for breeding offspring, and both parents are chosen independently and at random form this selection.
Then selections obtained via $k$-Pareto optimality-based sorting yield most offspring with parents offering choice.
\textit{Choice} here means that every parent is strictly superior to the other with respect to at least one objective, or both have the same values of all objectives.

\subsection{Kendall's $\tau$ Rank Correlation Coefficient and Statistical Tests}

Let us again consider the case of 2 continuous random variables introduced in \cref{sec:theoretical-exploration:random-vectors}.
Let us assume that $X_2=f(X_1)$ for some increasing function $f$. 
For almost all $(x, y)$, either $x R y$ or $y R x$ holds. Thus
$\divers(\Omega)=0$. Moreover, diversity only depends on the copula which encodes the dependency structure between $X_1$ and $X_2$, see \cref{sec:theoretical-exploration:random-vectors:marginal-independ,sec:theoretical-exploration:diversity}. 
Therefore, a value of
$\divers(\Omega)$ close to zero indicates $X_1$ and $X_2$ are strongly correlated via an increasing function. 
It leads to the idea that Kendall's $\tau$ rank correlation coefficient \cite{kendall1948rank} and diversity are strongly related concepts.

Let us consider a sample 
of $n$ points $X=\{(x_{i1}, x_{i2})_{i \in\{1, 2, ... , n\}}\}$. Duplicates  almost never occur and
the order in which points are drawn has no importance.
Therefore, $X$ should
be treated like a set. We  consider the counting measure $\#$
and the relation $R_{\llcorner}$ define on $X$. 
Diversity and choice of $X$ are denoted by
$\divers_{\#}$ and $\cho_{\#}$.

Kendall's $\tau$ correlation coefficient is defined as follows
$$\tau=\frac{\#\conco-\#\disco}{n_0},$$
where $\#\conco$ is the number of concordant pairs (pairs that do not offer choice in our terminology), $\#\disco$ is the number of discordant pairs (pairs that do offer choice), and $n_0$ is the total number of pairs. As duplicates are discarded and the pairs are not ordered, $n_0 = n(n-1)/2$.

From the above remarks, we have $\cho_{\#}=2 \# \disco + n$ and
$\#\conco+\#\disco = n_0$.
Combining these equalities, we obtain the following relation
$$\tau=\frac{n^2+n-2\cho_{\#}}{n^2-n}.$$
Dividing the numerator and the denominator by $n^2$ and then neglecting $\frac{1}{n}$, we obtains that approximately
$$\tau \approx 1-2\divers_{\#}.$$

This result means that the theory developed in this paper can be used for constructing non-parametric statistical tests generalizing Kendall's $\tau$ rank correlation coefficient and can be used for testing \textit{partial correlation}. 
Below, we illustrate this property by building an indicator for distinguishing between wealthy and non-wealthy states.
A group of states might be considered \textit{wealthy} if the following two conditions hold.
\begin{enumerate}
\item In the group there is no positive correlation between per capita income and the indicator representing education and health.   
\item If a state belongs to a group of wealthy states, then all states having higher per capita income and better value of education and health indicator, must also belong to that group.
\end{enumerate}

Now, we
define the set of all wealthy states as the largest group of states that are wealthy.
If Kendall's $\tau$ is used to compute correlation,
and correlation is considered to be positive if $\divers_{\#} < \frac{1}{2}$, then
the elements of the above set can be easily identified.

\begin{figure}[!h]
	\centering{
\begin{tikzpicture}

\begin{axis}[
legend cell align={left},
legend style={fill opacity=0.8, draw opacity=1, text opacity=1, draw=white!80!black,
at={(0.6,0.2)}, anchor=south west,
},
tick align=outside,
tick pos=left,
x grid style={white!69.0196078431373!black},
xlabel={Gross National Income (GNI), \$},
xmajorgrids,
xmin=-5878.905, xmax=136382.005,
xminorgrids,
xtick style={color=black},
xtick={0, 20000, 40000, 60000, 80000, 100000, 120000},
xticklabels={0, 20\,000, 40\,000, 60\,000, 80\,000, 100\,000, 120\,000},
scaled x ticks = false,
y grid style={white!69.0196078431373!black},
ylabel={Education  and  Health  Indicator},
ymajorgrids,
ymin=0.32215, ymax=0.97885,
ytick={.4, .5, .6, .7, .8, .9},
yminorgrids,
ytick style={color=black},
height=7cm,
width=13cm,
]
\addplot [only marks, mark=*, draw=black, fill=black, colormap/viridis]
table{%
x                      y
78162.3 0.925
67614.4 0.949
56364 0.939
42822.2 0.939
44518.9 0.925
129915.6 0.856
37065.2 0.921
32870.1 0.915
75064.8 0.912
46325.6 0.924
76075.2 0.8
44999.6 0.926
54264.9 0.917
72843.5 0.865
53245.1 0.92
43798 0.923
46250.8 0.913
62470.6 0.898
66203.3 0.84
42581.9 0.92
47979.5 0.858
51320.1 0.847
37930.8 0.91
37268 0.903
43608.8 0.893
31214.7 0.899
38085.4 0.897
34540.6 0.901
41243.3 0.896
38868.1 0.895
28664.2 0.89
33573 0.887
28144 0.878
32778.5 0.884
37236.4 0.824
24807.5 0.866
29458.5 0.856
26361.9 0.865
29499.6 0.856
34402.3 0.796
24117 0.855
26763.6 0.845
21665.5 0.847
};
\addlegendentry{wealthy}
\addplot [only marks, mark=*, draw=white!50.1960784313725!black, fill=white!50.1960784313725!black, colormap/viridis]
table{%
x                      y
26006.5 0.848
26103.6 0.843
23394.3 0.836
7455.1 0.775
28049.1 0.78
20291.3 0.827
22589.2 0.83
23286.1 0.804
20945.1 0.827
15409.6 0.807
24619.7 0.789
16261.3 0.794
19427.6 0.802
15629.4 0.796
23886.1 0.782
8855.8 0.769
14951.6 0.795
22093.1 0.794
13770.8 0.788
21565.1 0.792
19148.1 0.795
12202.1 0.776
20907.2 0.786
19470.2 0.788
22436.5 0.765
10788.9 0.766
10252.5 0.764
14006.4 0.776
17948 0.781
21517.2 0.592
16395.2 0.774
7361 0.743
16413.3 0.759
8189.1 0.743
5284.4 0.721
18704.5 0.767
16383.1 0.762
13311.6 0.763
15128.6 0.767
19043.6 0.697
10090.8 0.75
8244.8 0.736
16018 0.725
14145.1 0.754
11502.2 0.754
12405.2 0.748
14518.6 0.74
10111.3 0.742
13533 0.745
14662.8 0.698
11294.8 0.74
14303.1 0.716
10536.2 0.739
8350 0.73
13345.5 0.738
9791.4 0.735
14025.8 0.692
12762.2 0.727
10449.2 0.735
5371.5 0.704
5026.3 0.699
12756.4 0.722
5748.2 0.701
10095.8 0.726
10372.4 0.722
10248.6 0.725
3096.8 0.664
12087.2 0.666
11607.7 0.649
10382.7 0.701
7374.8 0.706
10064.1 0.691
10053.3 0.689
8395.1 0.682
8181.7 0.693
5256.5 0.684
5334.5 0.683
6154.9 0.674
7732 0.68
9769.8 0.64
2600.6 0.627
3291.5 0.638
6884.5 0.638
7194.9 0.647
7521.8 0.541
6049.1 0.648
7063.2 0.64
4746.7 0.645
7081.5 0.607
5663.5 0.624
4466.3 0.625
6290.9 0.533
5502.6 0.592
5371.2 0.606
2804.7 0.597
5442.9 0.527
2475.2 0.588
5048.6 0.586
3341.5 0.579
3463.7 0.579
5031.2 0.55
3838.7 0.579
4943.1 0.556
3069.8 0.574
3095.3 0.563
3527.3 0.513
2337.1 0.558
2880.7 0.555
3846.4 0.49
1319.7 0.512
2894.3 0.518
3318.8 0.497
2466.9 0.531
2440.6 0.536
1587.7 0.516
2711.8 0.516
1560.9 0.515
3216.2 0.473
3162.7 0.474
1616.9 0.498
1656.9 0.493
1670.2 0.493
1335 0.498
1262.2 0.487
1073.3 0.476
2250.1 0.494
1979.3 0.485
679.6 0.435
1870.8 0.479
2300.3 0.482
682.6 0.427
2218.4 0.442
1990.9 0.396
1882.5 0.418
1540.8 0.452
1523 0.448
1489.9 0.42
1529.2 0.42
691.2 0.404
1537.5 0.402
1369.5 0.424
1098.4 0.418
1058.2 0.414
889.5 0.353
587.5 0.352
};
\addlegendentry{non-wealthy}
\end{axis}

\end{tikzpicture}
		}
    \caption{Separation between wealthy and non-wealthy states based on $\divers_{\#}$.}
    \label{fig:rich-poor}
\end{figure}
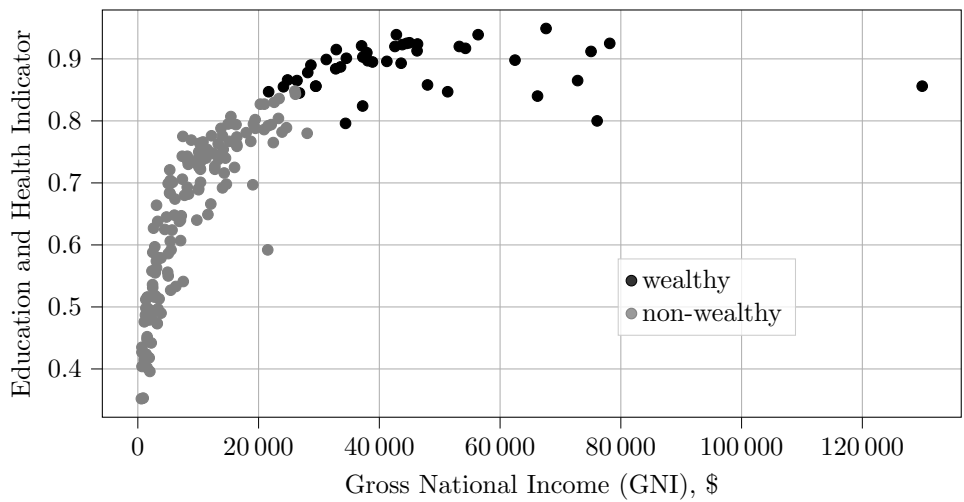

Indeed, Corollary~\ref{coro:max-div} says that the set of wealthy states must be a set of at-least $k$-Pareto optimal states for the relation \textit{higher income} and \textit{better education and health indicator}. 
For the year 2015\footnote{
Data source: United Nations Development Programme -
Human Development Reports \url{http://hdr.undp.org/en/data}.
},
we took Gross National Income (GNI) 
per capita at purchasing power parity (PPP) as the \textit{income indicator},
and the square root of the education and life expectancy as the \textit{education and health indicator}\footnote{
As suggested by the human development index, see \url{http://hdr.undp.org/sites/default/files/hdr2020_technical_notes.pdf}.
}. 
The scatter plot in \cref{fig:rich-poor} shows the resulting division of states into wealthy and non-wealthy.
We can observe that the wealthy states are defined as the states with GNI $\geq 20\,000\$$.
This seems perfectly plausible.
 
\subsection{Application 
to Recommender Systems}
\label{sec:applications:recommender-systems}

Let us again consider the housing example introduced in \cref{sec:pareto-optimality:definitions}.
If we aim to provide to the user a full set of possible alternative houses that might fit his preferences, then, according to Theorem~\ref{theorem:po-solution}, sorting available habitations by the increasing value of $\po$ is the best strategy.
As it was discussed in \cref{sec:practical-exploration:computation-complexity}, in the case of independent components of the underlying composite relations, the computation of $\po$ can be simplified by using tools from probability theory.
Apart from computational efficience, estimating $\po$ in this way 
has several additional advantages. 
\begin{itemize}
\item 
\textit{Such sorting results in fewer ties and a meaningful score.} 
Indeed, sorting items $x$ by increasing values of $\po(x)$, is the same as sorting by  decreasing values of $-\log(\po(x))$.
The self-information $-\log(F_i(x))$ \cite{jones1979elementary}, which is additive, indicates how much a characteristic $i$ is valued. 
In this case, there is no need to introduce any arbitrary coefficients as it is done when sorting by a weighted mean of the characteristics $x_i$.
\item 
If the condition of independence
holds, then \textit{rarer characteristics get valued more}. 
This makes sense from the economic point of view and is intuitively necessary for maximizing choice.
\item 
If beyond the relation $R$, there is complete uncertainty about the user's complex needs, tastes and desires, then
offering him a selection of maximum choice \textit{maximizes the likelihood he finds an appropriate item}.
\end{itemize}



\subsection{Constrained Multi-Objective Genetic Algorithms}

\subsubsection{Problem Definition}
A Multiobjective Constrained Optimization Problem (CMOP) is a mathematical problem that is defined as follows  \cite{Kumara2021RevisedGuidelines}: 

Minimize 
$$f_1(x), f_2(x), \dots, f_M(x)$$
subject to
\begin{align*}
    &g_i(x)\leq 0, i \in \{1, 2, ..., ng\},\\
    &h_j =0,  j \in \{ng+1, ng+2, ..., ng+nh\},\\
    &L_k \leq x_k \leq U_k, h \in \{1, \dots, D\},
\end{align*}
where 
\begin{itemize}
    \item $f_i$ represents the $i$-th objective function, 
    \item $M$ is the total number of conflicting objective functions,
    \item $x=(x_1, x_2, \dots, x_D)$ is a solution vector of length $D$,
    \item $L_k$ and $U_k$ are the lower and upper bounds of the search space at the $k$-th dimension.
\end{itemize}
Numerically, we consider a constraint  $h_j$ to be verified iff $h_j \in [-\epsilon, \epsilon]$.  
A solution is feasible iff all  $ng+nh$  constraints $g_i$ and $h_j$ are verified. 

\subsubsection{Problem Re-Definition with Preorder Relations}

In a more general setting, we can
represent a constraint $g_i \leq 0$ 
by the preorder relation
$R_{g_i}$ defined as follows: 
\begin{equation}
    x R_{g_i} y  \text{ iff }
    \begin{cases*}
      g_i(x) \leq 0 \text{ or }\\
      g_i(x) \leq g_i(y).
    \end{cases*}
\end{equation}

And a constraint $h_j  \in [a_j, b_j]$ 
can be represented by the preorder relation
$R_{h_j}$ defined as follows: 
\begin{equation}
    x R_{h_j} y  \text{ iff }
    \begin{cases*}
      h_j(x) \in [a_j, b_j]  \text{ or }\\
      h_j(y) \leq  h_j(x) \leq a_j \text{ or } \\
      b_j \leq h_j(x) \leq  h_j(y).
    \end{cases*}
\end{equation}

Then, the combination of
the constraints $g_i \leq 0, i \in \{1, 2, ..., ng\}$ and 
$h_j =0,  j \in \{ng+1, ng+2, ..., ng+nh\}$
can be represented by the  
preorder relation $R_c$ defined as follows.
$$x R_c y \text{ iff } x R_{g_1} y \text { and } \dots \text{ and } x R_{g_{ng}} y  
\text{ and } x R_{h_{ng+1}} y   \text{ and } \dots \text{ and } x R_{h_{ng+nh}} y.$$

The objective consisting in minimizing $f_i$ is represented by the preodrer relation $R_{f_i}$:
$$x R_{f_i} y \text{ iff } f_i(x) \leq f_i(y).$$
Minimization of all $M$ objectives $f_1, \dots, f_M$ is  represented by the preorder relation $R_f$ defined as follows
$$x R_f y  \text{ iff } x R_{f_1} y \text { and } \dots \text{ and } x R_{f_M} y.$$

Thus, the above CMOP can be represented by the lexicographic preorder relation
$$x R_{cf} y \text{\quad iff \quad}  x R^{*}_c y\text{ or }(x R^{=}_c y\text{ and  } x R_f y),$$
where  $x R^{*}_c y$ means ``$x R_c y$ and not $y R_c x$'', and $x R^{=}_c y$ means ``$x R_c y$ and $y R_c x$''. 
For the given $R_{cf}$,
constrained Pareto optimal solutions \cite{Kumara2021RevisedGuidelines} are solutions that are not 
Pareto dominated by any other solution. 

\subsubsection{Solution}

To solve the problem defined above,
we can use the standard Adaptive Differential Evolution Algorithm jDE \cite{jDEalgo} with $k$-Pareto optimality for  $R_{cf}$ as a fitness function. 
For any point $x$,
$k$-Pareto optimality of $x$ is the likelihood a point drawn at random from the population strictly Pareto dominates $x$ for $R_{cf}$. 
Smaller values of $\po$ mean better fitness. 
Under the independence assumption of  
objectives and constraints,
we can easily compute $k$-Pareto optimality $\po(x)$.
When saying $P(\{y |y R_{cf}^= x\})=0$, we assume the considered objectives and constraints are not constant on too large sets. 
Without this simplification, the  
computation becomes longer, see the derivation below. 
\begin{align*}
    \po(x) &= P(\{y |y R_{cf}^* x\}),\\
           &= P(\{y |y R_{cf} x\} - P(\{y |y R_{cf}^= x\}),\\
           &= P(\{y |y R_{cf} x\}-0,\\
           &= P(\{y |y R^{*}_c x\} \cup \{y |y R^{=}_c x\text{ and  } y R_f x\}),\\
           &= P(\{y |y R^{*}_c x\})+P(\{y |y R^{=}_c x\text{ and } y R_f x\}).\\
\end{align*}
Thus, if $x$ satisfies all constraints, which means $x R^{=}_c (0, \dots, 0, a_{ng+1}, \dots, a_{ng+nh})$, then
\begin{align*}
     \po(x) &= P(\{y |y R^{=}_c x\text{ and  } y R_f x\}),\\
            &= P(\{y |y R^{=}_c (0, \dots, 0, a_{ng+1}, \dots, a_{ng+nh})\})P(\{y |y R_f x\}).
\end{align*}

Now, let $F_i(z)=P(\{y | f_i(y)\leq z\})$ be the cumulative probability distribution of~$f_i$.
Then,
\begin{align*}
P(\{y |y R_f x\})&=P\left(\bigcap_{i \in \{1, \dots, M\}}\{y |f_i(y) \leq f_i(x)\}\right) ,\\ 
&=\prod_{i \in \{1,\dots,M\}}P(y | f_i(y) \leq f_i(x)).\\
&=\prod_{i \in \{1,\dots,M\}}F_i(f_i(x)).
\end{align*}
Otherwise, if at least one constraint is not satisfied by $x$, then $P(\{y |y R^{=}_c x\})=0$ and
\begin{align*}
      \po(x)&=P(\{y |y R^{*}_c x\}) \\
            &=P(\{y |y R_c x\}) \\ 
            &=\prod_{i \in \{1,\dots,ng\}} G_i(g(x)) 
              \prod_{j \in \{ng+1,\dots,ng+nh\}}P(y |y R_{h_j} x). \label{algo:1}
\end{align*}

Any cumulative probability distribution defined above, $F_i$ for $f_i$ and $G_i$ for $g_i$, 
can be estimated via its empirical cumulative probability distribution. 
Note, for a population $\{x_1, \dots, x_k, \dots ,x_{ps}\}$ of size $ps$, and for any real valued function $f$, the empirical cumulative probability distribution $\hat{F}$ of $f$ is defined as follows:
$$\hat{F}(z)=\frac{\#(\{x_k| f(x_k) \leq z\})}{ps}.$$
In a similar way, $P(y |y R_{h_j} x)$ can be estimated by $\widehat{H_j}^*(h_j(x))$, where  $\widehat{H_j}^*(z)$ is defined as
\begin{equation}
\label{eq:cumulative-H}
    \widehat{H_j}^*(z)=
    \begin{cases*}
      \frac{\#(\{x_k| z \leq h_j(x_k) \leq b_j\}}{ps}  \text{ \quad if $z<a_j$},\\
      \frac{\#(\{x_k|  a_j \leq h_j(x_k) \leq b_j]\}}{ps}  \text{\quad if $a_j \leq z \leq b_j$},\\
      \frac{\#(\{x_k| a_j \leq h_j(x_k) \leq z\}}{ps}  \text{ \quad if $z>b_j$}.
    \end{cases*}
\end{equation}

The computation of every $\hat{F_i}$ can be performed as follows: 
\begin{itemize}
    \item sort the values $f(x_k)$ in increasing order, and store them in an array;
    \item create two new arrays;
    \item loop over the sorted values $f(x_k)$; each time a new distinct value $f(x_k)$ is encountered:
    \begin{itemize}
    \item append the previously encountered $f(x_{k})$ to the first array, 
    \item append to the second array the loop counter, which is equal to the value of the empirical cumulative probability distribution $\hat{F_i}$ of the previously encountered value $f(x_{k})$.  
    \end{itemize}
\end{itemize}

Thus, retrieving  $\hat{F_i}(x)$ can be performed via a binary lookup with run time $O( \log ps)$.
Computation of $\widehat{H_j}^*$ can be performed in the same way.
In this case, all three cases of the definition in \cref{eq:cumulative-H} are treated separately. 
Moreover, we have the estimation
\begin{equation*}
P(\{y |y R^{=}_c (0, \dots, 0, a_{ng+1}, \dots, a_{ng+nh})\})
=\prod_{i \in \{1,\dots,ng\}}\hat{G}_i(0) \prod_{j \in \{ng+1,\dots,ng+nh\}} \widehat{H_j}^*(a_j).
\end{equation*}

Finally, 
it is possible to show that
the total run time of the $k$-Pareto optimality based sorting is $O((ng+nh+M) ps \log ps)$.

\subsection{Exploratory Database Queries}



Simple database queries $q$, objectives, and constraints in optimization problems often consist in requiring a continuous attribute to be in a given interval, or a discrete attribute to be equal to a given value. 
Conceptually, those queries are boolean functions. 
Complex queries are often conjunctions of the form $r=q_1 \land q_2 \land \dots \land q_n$. 

In our formalism, these simple queries translate 
into simple pre-order relations of the form $x R y$. 
Requiring an element to be in an interval can be represented by $x R_q y$ iff $x$ is in the desired interval, or $x$ is not situated further from the interval than $y$\footnote{
Hence, the strict version of this relation is defined as follows: $x R_q^* y$ iff $x$ is in the desired interval and $y$ is not, or if $x$ is situated closer to the desired interval.
}. 
Requiring an attribute to be equal to a given value translates into the relation $x R_{q'} y$ if for $x$ the attribute takes the required value\footnote{The strict version of this relation is defined as $x R_{q'}^* y$ if for $x$ the attribute takes the required value but not for $y$.}. 
Complex queries then translate into the composite relations of the form 
$R_r=R_{q_1} \land R_{q_2} \land \dots \land R_{q_n}$. The simple sub-relations $R_{q_2}$ 
are pre-order relations, and, therefore, $R_r$ is also a pre-order relation and is transitive. 
However, these relations are not partial order relations, as reflexivity does not necessarily hold. 
A "topological" sorting
according to our partial order relation $R_r$ can be viewed as a valid fuzzy relaxation of the strict functional query. 
There are many possible fuzzy relaxations and the problem is to find one that is suitable for a given application. 
The $k$-Pareto optimality is one of such fuzzy extensions of the query.
It is $0$ if all criteria are satisfied, and higher values of $k$-Pareto optimality indicate worse results. 
The maximum choice theorem applies here, and the user is offered the maximum choice.
This is of particular interest for exploratory queries, such as job search,
especially, if there are no items in the database that satisfy all the criteria. 
Direct brute-force search for selections offering maximum choice is unfeasible as there are too many selections to consider. 

Moreover, in the above formalism, one can treat classical optimization objectives in the same way. 
Maximizing an attribute $x$ can be represented by the relation $xRy$ iff $x \geq y$, and the minimization can be represented by the relation $\leq$.  
In the above framework, negation can be represented via the relation $R^{-1}$ defined 
by $x R^{-1} y$ iff $y R x$.

\subsection{Scheduling Algorithms}

In the case of scheduling algorithms, $xRy$ can  be  given  the meaning  
`$x$ depends on $y$'. 
Then, selections represent sets of tasks that remain to be processed. 
Having a large choice means \textit{having much freedom to parallelize tasks} or \textit{having flexibility in case the rescheduling is required}.



\end{document}


%

%

\onecolumn
\aistatstitle{Instructions for Paper Submissions to AISTATS 2022: \\
Supplementary Materials}

\section{FORMATTING INSTRUCTIONS}

To prepare a supplementary pdf file, we ask the authors to use \texttt{aistats2022.sty} as a style file and to follow the same formatting instructions as in the main paper.
The only difference is that the supplementary material must be in a \emph{single-column} format.
You can use \texttt{supplement.tex} in our starter pack as a starting point, or append the supplementary content to the main paper and split the final PDF into two separate files.

Note that reviewers are under no obligation to examine your supplementary material.

\section{MISSING PROOFS}

The supplementary materials may contain detailed proofs of the results that are missing in the main paper.

\subsection{Proof of Lemma 3}

\textit{In this section, we present the detailed proof of Lemma 3 and then [ ... ]}

\section{ADDITIONAL EXPERIMENTS}

If you have additional experimental results, you may include them in the supplementary materials.

\subsection{The Effect of Regularization Parameter}

\textit{Our algorithm depends on the regularization parameter $\lambda$. Figure 1 below illustrates the effect of this parameter on the performance of our algorithm. As we can see, [ ... ]}

\vfill